\documentstyle[11pt]{article}
   
        \newcounter{sect}

        \newcommand{\be}{\begin{equation}}

        \newcommand{\ee}{\end{equation}}

        \newcommand{\bea}{\begin{eqnarray}}

        \newcommand{\eea}{\end{eqnarray}}

        \newcommand{\nno}{\nonumber \\}

        \newtheorem{theorem}{Theorem}

        \newtheorem{definition}{Definition}

        \newcommand{\A}{{\cal A}}

        \newcommand{\D}{{\cal D}}\newcommand{\HH}{{\cal H}}

        \newcommand{\LL}{{\cal L}}

        \newcommand{\B}{{\cal B}}

        \newcommand{\K}{{\cal K}}

        \newcommand{\s}{\sigma}

        \newcommand{\p}{\partial}

        \newcommand{\dd}{|\D|}

        \newcommand{\n}{\parallel}  

\newcommand{\bma}{\left(\begin{array}{cc}}

\newcommand{\ema}{\end{array}\right)}

\newcommand{\bca}{\left(\begin{array}{c}}

\newcommand{\eca}{\end{array}\right)}

\newcommand{\sr}{\stackrel}

\newcommand{\la}{\langle}

\newcommand{\ra}{\rangle}

        \newcommand{\R}{\mathbf R}

        \newcommand{\C}{\mathbf C}

\newcommand{\Z}{\mathbf Z}

\newcommand{\NN}{{\mathcal N}}

\newcommand{\ben}{\begin{displaymath}}

        \newcommand{\een}{\end{displaymath}}

        \newcommand{\bean}{\begin{eqnarray*}}

        \newcommand{\eean}{\end{eqnarray*}}

\newcommand{\proof}[1]{\noindent \normalfont{\textbf{Proof} \  \  #1 \hfill 
$\Box$} \par}

     \newtheorem{lemma}[theorem]{Lemma}

      \newtheorem{proposition}[theorem]{Proposition}

        \newtheorem{corollary}[theorem]{Corollary}

\begin{document}
\titlepage

\begin{center}{\bf THE HOCHSCHILD CLASS OF THE CHERN CHARACTER FOR SEMIFINITE SPECTRAL TRIPLES}
\\
\vspace{.5 in}

\title{}
\author{}
{\bf Alan L. Carey}\\Mathematical Sciences Institute\\
Australian National University\\
Canberra, ACT. 0200, AUSTRALIA\\
e-mail: acarey@maths.anu.edu.au\\
\vspace{.2 in}

{\bf John Phillips}\\Department of Mathematics and Statistics\\
University of Victoria\\Victoria, B.C. V8W 3P4, CANADA\footnote{Address for
correspondence}\\
e-mail: phillips@math.uvic.ca\\
\vspace{.2 in}

{\bf Adam Rennie}\\
School of Mathematical and Physical Sciences\\
University of Newcastle\\
Callaghan, NSW, 2308 AUSTRALIA\\
e-mail: adam.rennie@newcastle.edu.au\\

\vspace{.2 in}
{\bf and}
\vspace{.2 in}

{\bf Fyodor A. Sukochev}\\
School of Informatics and Engineering\\
Flinders University\\
Bedford Park S.A 5042 AUSTRALIA\\
e-mail: sukochev@maths.flinders.edu.au\\
\vspace{.25 in}

RUNNING TITLE: THE HOCHSCHILD CLASS OF THE CHERN CHARACTER 

All authors were supported by grants from ARC (Australia) and
NSERC (Canada), in addition the third named author acknowledges a 
University of Newcastle early career researcher grant.

\end{center}

\newpage
\begin{abstract}

We provide a proof of Connes' formula for a representative of the Hochschild 
class of the Chern character for $(p,\infty)$-summable spectral triples. 
Our proof is valid for all semifinite von Neumann algebras, and all 
integral $p\geq 1$. We employ the minimum possible hypotheses on the spectral 
triples.
\footnote{
von Neumann algebra, Fredholm module, cyclic cohomology, chern character.}
\end{abstract}
\newpage
\section{Introduction}

A key result in the quantised calculus of Alain Connes 
(\cite[IV.2.$\gamma$]{BRB})
is the formula for the Hochschild class of the Chern 
character of a $(p,\infty)$-summable spectral triple (these notions are 
explained below).

Our aim is to generalise this formula to encompass the situation in which
one uses, instead of the bounded operators on  Hilbert space and its 
various ideals
of compact operators, a general semifinite von Neumann algebra
and the analogous ideals as described for example in \cite{FK,S}.
Moreover we aim to prove the formula in the 
greatest possible generality with only the absolutely essential side
conditions. This is a delicate matter as regards the amount of smoothness or 
regularity necessary. The result has been stated, \cite[Theorem 10.32]{Va}, 
with the hypothesis that the algebra be `twice quantum differentiable' 
(see below), but the proof appearing in \cite[pp 470-479]{Va} does not quite 
hold with this hypothesis. We employ the hypothesis of `$\max\{2,p-2\}$ quantum 
differentiability', and while this is sufficient, the necessity of this 
condition is unknown to us. Indeed, we only require this stronger hypothesis 
at one (crucial) point, Proposition \ref{usesomega}, but we isolate the 
particular statement which uses this hypothesis in Lemma \ref{bounded}.

A rationale for the extension of Connes'
spectral geometry to the case of general 
semifinite von Neumann algebras is presented in \cite{BF}.
Examples where this notion arises naturally
are non-smooth foliations \cite{BF,P}, the
$L^2$-index theorem (see \cite{Ma} and references therein)
and $L^2$ spectral flow \cite{CP1,CP2}.

In order to describe our results some preliminary machinery
is needed (all of this is contained in \cite{BRB} in the type I case). We deal 
with this in Section \ref{BM}. We first describe spectral triples for 
semifinite von Neumann algebras, including definitions of smoothness and 
summability. We then briefly recall the Hochschild and cyclic cohomology 
theories, and explain what the Hochschild class of the Chern character is, 
and what kind of information it contains.

The last preliminary subsection describes results from \cite{CPS2}, where a 
proof of the connection between the trace of the heat kernel and the Dixmier 
trace is presented. The idea has previously appeared, \cite[p 563]{BRB}, 
but this is the first proof simultaneously valid for the case $p=1$ and the 
general semifinite case. It is a key tool in our proof.

Section \ref{MR} begins with a statement of the main result and its main 
corollaries. The expert reader can skip straight to Subsections \ref{31} 
and \ref{32} for our result, its corollaries, and how it relates to the 
significant body of previous work on this general topic. We then set out 
the proof as clearly as possible, and in the greatest possible generality.  
The proof is considerably simplified by the assumption of invertibility of 
the `Dirac' operator $\D$, but a standard construction in $K$-homology and 
computations contained in the Appendix show that the result remains true 
even when this is not the case.

\section{Background Material and Preliminary Results}\label{BM}

\subsection{Spectral Triples}

We begin with some semifinite versions of standard definitions and results.

\begin{definition} A semifinite
spectral triple $(\A,\HH,\D)$ is given by a Hilbert space $\HH$, a 
$*$-algebra $\A\subset \NN$ where $\NN$ is a semifinite von Neumann algebra 
acting on
$\HH$, and a densely defined unbounded self-adjoint operator $\D$ affiliated 
to $\NN$ such that

1) $[\D,a]$ is densely defined and extends to a bounded operator in
$\NN$ for all $a\in\A$

2) $(\lambda-\D)^{-1}\in\K(\NN)$ for all $\lambda\not\in{\R}$

Here $\K(\NN)$ is the ideal of $\tau$-compact operators in $\NN$ (this is 
explained in the next section). 
We say that $(\A,\HH,\D)$ is even if in addition there is a ${\Z}_2$-grading 
such that $\A$ is even and $\D$ is odd. That is an operator $\Gamma$ such that 
$\Gamma=\Gamma^*$, $\Gamma^2=1$, $\Gamma a=a\Gamma$ for all $a\in\A$ and 
$\D\Gamma+\Gamma\D=0$. Otherwise we say that $(\A,\HH,\D)$ is odd.
\end{definition}

{\bf Remark} We will write $\Gamma$ in all our formulae, with the 
understanding 
that if 
$(\A,\HH,\D)$ is odd, $\Gamma=1$ and of course, we drop the assumption that 
$\D\Gamma+\Gamma\D=0$. Alas, we will also employ the gamma function  
in this paper, but the meaning of the symbol `$\Gamma$' should be clear from 
context. Henceforth we omit the term semifinite
as it is implied by the use of a faithful normal semifinite trace $\tau$
on $\NN$ in all of the subsequent text.

\begin{definition} A spectral triple $(\A,\HH,\D)$ is $QC^k$ for $k\geq 1$ 
($Q$ for quantum) if for all $a\in\A$ 
the operators $a$ and $[\D,a]$ are in the domain of $\delta^k$ where 
$\delta(T)=[\dd,T]$ is the (partially defined) derivation 
on $\NN$ defined by $\dd$. We say that 
$(\A,\HH,\D)$ is smooth if it is $QC^k$ for all $k\geq 1$.
\end{definition}

{\bf Remark} The notation is meant to be analogous to the classical case, but 
we introduce 
the $Q$ so that there is no confusion between quantum differentiability of 
$a\in\A$ and classical differentiability of functions. We may also speak about 
a $QC^0$ spectral triple, where only $a$ and $[\D,a]$ are assumed bounded. We 
also note that if $T\in{\mathcal N}$, one can show that $[\dd,T]$ is bounded 
if and only if $[(1+\D^2)^{1/2},T]$ is bounded, by using the functional 
calculus to show that $\dd-(1+\D^2)^{1/2}$ is a bounded operator and
lies in $\NN$.

\subsubsection{Summability}

Recall from \cite{FK} that if $S\in\mathcal N$ the {\bf t-th generalized
singular value} of $S$ for each real $t>0$ is given by
$$\mu_t(S)=\inf\{\n SE\n\ \vert \ E \mbox{ is a projection in }
{\mathcal N} \mbox { with } \tau(1-E)\leq t\}.$$
We write
$T_1\prec\prec T_2$ to mean that $\int_0^t \mu_s(T_1)ds\leq
\int_0^t \mu_s(T_2)ds$ for all $t>0$.

\begin{definition}\label{opideal}
If $\mathcal I$ is a $*$-ideal in $\mathcal N$
which is complete in a norm $\n\cdot\n_{\mathcal I}$
then we will call $\mathcal I$ a {\bf symmetric
operator ideal} if\\
(1) $\n S\n_{\mathcal I}\geq \n S\n$ for all $S\in \mathcal I$,\\
(2) $\n S^*\n_{\mathcal I} = \n S\n_{\mathcal I}$ for all $S\in \mathcal I$,\\
(3) $\n ASB\n_{\mathcal I}\leq \n A\n \:\n S\n_{\mathcal I}\n B\n$ for all
$S\in \mathcal I$, $A,B\in \mathcal N$.\\
Since $\mathcal I$ is an ideal in a von Neumann algebra, it follows
from I.1.6, Proposition 10 of \cite{Dix} that if $0\leq S \leq T$
and $T \in {\mathcal I}$, then $S \in {\mathcal I}$ and $\n S\n_{\mathcal I}
\leq \n T\n_{\mathcal I}$. 
\end{definition}

Such ideals are special cases of symmetric operator spaces (see \cite{S} and 
references therein).
The main examples of such ideals that we consider in this paper
are the spaces
$${\mathcal L}^{(1,\infty)}({\mathcal N})=
\left\{T\in{\mathcal N}\ |\ \Vert T\Vert_{_{{\mathcal L}^{(1,\infty)}}}
:=   \sup_{t> 0}
\frac{1}{\log(1+t)}\int_0^t\mu_s(T)ds<\infty\right\}.$$
and with $p>1$,
$$\psi_p(t)=\left\{\begin{array}{ll} t & \mbox{for } 0\leq t\leq 1\\
                                     t^{1-\frac{1}{p}} & \mbox{for } 1\leq t
\end{array}\right.$$
$${\mathcal L}^{(p,\infty)}({\mathcal N})
=\left\{T\in{\mathcal N}\ |\ \Vert T\Vert_{_{{\mathcal L}^{(p,\infty)}}}
:= \sup_{t> 0}
\frac{1}{\psi_p(t)}\int_0^t\mu_s(T)ds<\infty\right\}.$$
For $p>1$ there is also the equivalent definition (see for example 
\cite[Section 5]{S})
$${\mathcal L}^{(p,\infty)}({\mathcal N})
=\left\{T\in{\mathcal N}\ |\ \sup_{t> 0}
\frac{t}{\psi_p(t)}\mu_t(T)<\infty\right\}.$$

It is well-known (see e.g. \cite{S}) that for 
$T_1\in {\mathcal N}$, $T_2\in {\mathcal L}^{(p,\infty)}({\mathcal N})$, 
$p\in [1,\infty)$, 
the condition $T_1\prec\prec T_2$ implies that 
$T_1\in {\mathcal L}^{(p,\infty)}({\mathcal N})$. We denote the norm 
on $\LL^{(p,\infty)}$ by $\n\cdot\n_{(p,\infty)}$.

As we will not change ${\mathcal N}$ throughout the paper we will suppress
the $({\mathcal N})$ to lighten the notation. The reader should note
that ${\mathcal L}^{(p,\infty)}$ is often taken to mean an ideal in
the algebra $\widetilde{\mathcal N}$ of $\tau$-measurable operators affiliated
to ${\mathcal N}$. Our notation is however consistent with that of \cite{BRB}
in the special case ${\mathcal N}={\mathcal B}({\mathcal H})$. With this 
convention the ideal of $\tau$-compact operators, ${\mathcal
  K}({\mathcal N})$, 
consists of those $T\in{\mathcal N}$ (as opposed to $\widetilde{\mathcal N}$) 
such that 
\ben \mu_\infty(T):=\lim _{t\to \infty}\mu_t(T)  = 0.\een

\begin{definition}
A spectral triple $(\A,\HH,\D)$
is called $(p,\infty)$-summable if 
$(1+\D^2)^{-1/2}\in {\mathcal L}^{(p,\infty)}$.
\end{definition}

We will also require the ideals $\LL^p({\mathcal N})$ 
and $\LL^{(p,1)}({\mathcal N})$, for $p\geq 1$. An 
operator $T\in {\mathcal N}$ is in $\LL^p({\mathcal N})$ if
\ben \n T\n_p:=\tau( |T|^p)^{1/p}<\infty.\een
In the Type I setting these are the usual Schatten ideals. Again we will 
simply write $\LL^p$ for these ideals in order to simplify the notation, and 
denote the norm on $\LL^p$ by $\n\cdot\n_p$. An operator $T\in{\mathcal N}$ 
is in $\LL^{(p,1)}({\mathcal N})$ if, \cite{S},
\ben \n T\n_{(p,1)}:=(1/p \int_0^\infty (t^{1/p}\mu_t(T)) dt/t)  <\infty.\een
For $p=1$ the ideal $\LL^{(p,1)}$ coincides with $\LL^1$. We denote the norm 
on $\LL^{(p,1)}$ by $\n\cdot\n_{(p,1)}$. If $1<p<\infty$ 
and $\frac{1}{p}+\frac{1}{q}=1$, then the K\"{o}the dual of $\LL^{(p,1)}$ 
is $\LL^{(q,\infty)}$, \cite{DDP2}. For $p=1$ the K\"{o}the dual of $\LL^1$ 
is just ${\mathcal N}$.

We use the following results repeatedly. They tell us the summability of 
various operators associated to a $(p,\infty)$-summable spectral triple. The
results are established in \cite[Propositions 1.1 and 1.2]{CS}, 
namely that for any $\tau$-measurable operators $T_1$ and $T_2$ we have
\ben \mu(T_1T_2)\prec\prec\mu(T_1)\mu(T_2),\een
where $\mu(T)$ denotes the function $s\to\mu_s(T)$. Moreover, for any 
self-adjoint $\tau$-measurable operators $T$ and $S$ with $T\geq 0$, 
$$-T\leq S\leq T \ \ implies \ \ S\prec\prec T.$$

\begin{lemma}[\cite{Va}]\label{sizes} Let $(\A,\HH,\D)$ be 
a $(p,\infty)$-summable 
$QC^k$ spectral triple, $k\geq 0$, 
with $p\geq 1$ and $\D$ invertible. Then for all $a\in\A$
\ben [F,a],\ [F,\delta(a)],\ \ \cdots,\ [F,\delta^k(a)]\in 
\LL^{(p,\infty)},\een
where $F=\D\dd^{-1}$.

\proof{We start with the formula
\ben \dd^{-1}=\frac{1}{\pi}\int_0^\infty(\lambda+\D^2)^{-1}
\frac{d\lambda}{\sqrt{\lambda}}.\een
This is used to rewrite $[F,a]$ in the following way.
\bean [F,a]&=&[\D\dd^{-1},a]=[\D,a]\dd^{-1}+\D[\dd^{-1},a]\nno
&=& \frac{1}{\pi}\int_0^\infty\left([\D,a](\lambda+\D^2)^{-1}+
\D[(\lambda+\D^2)^{-1},a]\right)\frac{d\lambda}{\sqrt{\lambda}}\nno
&=&\frac{1}{\pi}\int_0^\infty\left([\D,a](\lambda+\D^2)^{-1}-\D^2
(\lambda+\D^2)^{-1}[\D,a](\lambda+\D^2)^{-1}\right.\nno
& & \qquad\qquad\qquad\qquad\qquad\left.-\D(\lambda+\D^2)^{-1}[\D,a]
\D(\lambda+\D^2)^{-1}\right)\frac{d\lambda}{\sqrt{\lambda}}\nno
&=&\frac{1}{\pi}\int_0^\infty\left(\lambda(\lambda+\D^2)^{-1}[\D,a]
(\lambda+\D^2)^{-1}-
\D(\lambda+\D^2)^{-1}[\D,a]\D(\lambda+\D^2)^{-1}\right)
\frac{d\lambda}{\sqrt{\lambda}}.
\eean

The second last equality comes from \cite[Lemma 2.3]{CP1}, whose proof 
requires only $QC^0$, as opposed to the usual resolvent calculation which 
requires $QC^1$. The final equality comes from 
$\D^2(\lambda+\D^2)^{-1}=1-\lambda(\lambda+\D^2)^{-1}$. We now 
suppose that $a^*=-a$ so that $[\D,a]^*=[\D,a]$ and similarly for $[F,a]$. 
Then we may 
employ the inequality $$-\n[\D,a]\n T^*T\leq T^*[\D,a]T\leq \n[\D,a]\n T^*T$$ 
for all 
$T\in\NN$. Applying this inequality under the above integral yields
\bean [F,a] &\leq &\frac{\n[\D,a]\n}{\pi}
\int_0^\infty(\lambda+\D^2)^{-1}\frac{d\lambda}
{\sqrt{\lambda}}\nno
&=&\n[\D,a]\n\dd^{-1},\eean
and similarly $[F,a]\geq -\n[\D,a]\n\dd^{-1}$. Thus 
$[F,a]\prec\prec\n[\D,a]\n\dd^{-1}$, in particular 
$[F,a]\in\LL^{(p,\infty)}$,
and by linearity this is true for all $a\in\A$. However, for not necessarily 
self-adjoint $a\in\A$, the precise inequality is
\be [F,a]\prec\prec (\n[\D,Re(a)]\n+\n[\D,Im(a)]\n)\dd^{-1}.
\label{precineq}\ee
From the comments in Definition \ref{opideal}, this shows that 
$[F,a]\in\LL^{(p,\infty)}$.
The remainder of the  result is proved using the same  argument by 
replacing $a$ by $\delta^i(a)$, for $i=1,...,k$, and using the boundedness 
of $[\D,\delta^i(a)]=\delta^i([\D,a])$.}
\end{lemma}

The following lemma is a consequence of the previous result. This is the  
point at which more smoothness than $QC^2$ is required. 
The analogous statement 
in \cite[Lemma 10.27]{Va}, is a little lax about the degree of 
smoothness necessary to perform the iterated commutators with $\dd$ in the 
proof.

\begin{lemma}\label{bounded} Let $p\geq 1$ and $k=\max\{1,p-2\}$. Suppose 
that $(\A,\HH,\D)$ is a $QC^k$ $(p,\infty)$-summable spectral triple with 
$\D$ invertible. For 
all $a_0,...,a_{p-1}\in\A$, and $T\in\NN\cap{\rm dom}(\delta)$, the operators
\ben \dd^{p-2}a_0[F,a_1]\cdots[F,a_{p-1}]
FT\dd\ \ \ {\rm and}\ \ \ \dd Ta_0[F,a_1]\cdots[F,a_{p-1}]F\dd^{p-2}\een
are densely defined and bounded (or, more accurately, extend to bounded 
operators).

\medskip
\proof{The proof is essentially the same as that in \cite{Va}. First, the 
triple is at least $QC^1$, so $[F,a]\dd=[\D,a]-F\delta(a)$ is bounded for 
all $a\in\A$. This allows one to check the cases $p=1,2$. For $p>2$ we have
\bean &&\dd^{p-2}a_0[F,a_1]\cdots[F,a_{p-1}]FT\dd \nno
&=&\sum_{j=0}^{p-2}\left(\begin{array}{c}p-2\\j\end{array}\right)
\delta^j(a_0)\dd^{p-2-j}[F,a_1]\cdots[F,a_{p-1}]F(\dd T-\delta(T)),\eean
and similarly for the other operator. Now $[F,a_{p-1}](\dd T-\delta(T))$ is 
bounded and in $\NN$, since the triple is $QC^1$ and
$T\in\mbox{dom}(\delta)$. 
Similarly
$j=p-1,p-2$, $\dd^{p-2-j}[F,a_1]$ is bounded and in $\NN$, and 
of course $\delta^j(a_0)$ is 
bounded for $0\leq j\leq p-2$ since the triple is $QC^{p-2}$. Hence we may 
consider only those terms with $0\leq j\leq p-3$.

One now continues to take commutators, observing that we obtain a sum of 
bounded operators in $\NN$ plus the term
\ben a_0\dd[F,a_1]\dd[F,a_2]\cdots\dd[F,a_{p-2}][F,a_{p-1}]\dd FT\een
which is bounded and in $\NN$ since the triple is $QC^1$. Similar comments apply to the 
second operator.}
\end{lemma}

\subsection{Hochschild and Cyclic Cohomology}\label{hochscyclic}

For a locally convex unital algebra $\A$, we denote by $C^n(\A)$ the linear 
space of continuous $n+1$-multilinear functionals on $\A^{n+1}$. The 
Hochschild coboundary of $\phi\in C^n(\A)$ is the functional $b\phi\in 
C^{n+1}(\A)$ defined by
\bean (b\phi)(a_0,...,a_{n+1})&=& \phi(a_0a_1,a_2,...,a_{n+1})\nno
&+& \sum_{i=1}^n(-1)^i\phi(a_0,...,a_ia_{i+1},...,a_{n+1})\nno
&+& (-1)^{n+1}\phi(a_{n+1}a_0,a_1,...,a_n),\ \ \ a_0,...,a_{n+1}\in\A.\eean
One can easily check that $b^2=0$. The Hochschild cohomology, denoted 
$HH^*(\A,\A^*)$, is then the cohomology of the complex $(C^*(\A),b)$. The 
notation is explained in \cite{L}. We denote the space of Hochschild 
$k$-cocycles by $Z^k(\A)$.

Let $C^n_\lambda(\A)$ be the subspace of $C^n(\A)$ consisting of 
functionals $\phi$ such that
\ben \phi(a_0,...,a_n)=(-1)^n\phi(a_n,a_0,...,a_{n-1}),\ \ \ a_0,...,a_n\in\A.
\een
The Hochschild coboundary maps $C^n_\lambda(\A)$ to $C^{n+1}_\lambda(\A)$, so 
we can define the cyclic cohomology of $\A$, denoted $HC^*(\A)$, to be the 
cohomology of the complex $(C^*_\lambda(\A),b)$. The space of cyclic cocycles 
is denoted $Z_\lambda(\A)$.

Connes shows that there is a long exact sequence, \cite[III.1.$\gamma$]{BRB},
\ben \cdots\sr{B}{\to}HC^p(\A)\sr{S}{\to}
HC^{p+2}(\A)\sr{I}{\to}HH^{p+2}(\A,\A^*)
\sr{B}{\to}HC^{p+1}(\A)\sr{S}{\to}\cdots\een
where $I$ is the map induced on cohomology by the inclusion of 
complexes $C^n_\lambda(\A)\hookrightarrow C^n(\A)$. The operator $B$ will 
not concern us in this paper, see \cite{L} and \cite[III.1.$\gamma$]{BRB}, 
however the periodicity operator $S$ is important for three reasons. The 
first reason is that cyclic cohomology groups are filtered by powers of $S$, 
so in general $HC^k(\A)$ consists of (classes of) sums 
\ben \phi_k+S\phi_{k-2}+S^2\phi_{k-4}+
...+S^{[k/2]}\phi_{k-2[k/2]}\een
where $\phi_j$ is cyclic, $b\phi_j=0$, $j=k-2[k/2],...,k$.
As $\mbox{Image}S=\ker I$, we have 
\ben I(\phi_k+...+S^{[k/2]}\phi_{k-2[k/2]})=I(\phi_k).\een

Consequently, pairing a cyclic cocycle $\phi$ with a Hochschild cycle yields 
the 
same result as pairing the Hochschild class of the cyclic cocycle $\phi$ with a
Hochschild cycle. Equivalently, any cocycle in the image of $S$ has zero 
Hochschild class.

Secondly, because $S:HC^n(\A)\to HC^{n+2}(\A)$ for all $n$, we may define even 
periodic cyclic cohomology as the inductive limit 
$H^{ev}(\A):=\lim_{\to}(HC^{2n}(\A),S)$, and similarly the odd periodic 
cyclic cohomology as $H^{odd}(\A):=\lim_{\to}(HC^{2n+1}(\A),S)$.

Thirdly, we use  constructions closely related to the periodicity operator in 
the Appendix to complete the proof of our main result for the case where our 
spectral triple $(\A,\HH,\D)$ has $\D$ noninvertible.

The (dual) homology theories require more care regarding the completeness of 
$\A$ and the appropriate tensor product. We only need Hochschild cycles (not 
their classes), so we may ignore these difficulties in this paper. We only 
require the definition of the Hochschild boundary. If 
$c=\sum_{i=1}^n a^i_0\otimes a^i_1\otimes\cdots\otimes a^i_k$ then
\bean bc=b(\sum_{i=1}^na^i_0\otimes a^i_1\otimes\cdots\otimes a^i_k)&=& 
\sum_{i=1}^n a^i_0a^i_1\otimes a^i_2\otimes\cdots\otimes a^i_k\nno
&+&\sum_{i=1}^n\sum_{j=1}^{k-1}(-1)^ja^i_0\otimes\cdots\otimes a^i_ja^i_{j+1}
\otimes\cdots\otimes a^i_k\nno
&+& (-1)^k\sum_{i=1}^na^i_ka^i_0\otimes\cdots\otimes a^i_{k-1}.\eean
We say that $c$ is a Hochschild cycle if $bc=0$. When the Hochschild homology 
is well-defined we denote it by $HH_*(\A)$.

An important point is that if $\phi$ is a  $k-1$-multilinear functional, and 
$c$ is a Hochschild $k$-cycle, $(b\phi)(c)=\phi(bc)=0$, so the pairing of a 
Hochschild coboundary with a Hochschild cycle vanishes. This follows 
immediately from the definitions.

Our only result in this section consists of a mild generalisation of a 
standard result for the 
behaviour of Hochschild (co)homology with respect to derivations. This 
result will simplify many later computations.

\begin{lemma}\label{hochs} Let $\NN$ be a semifinite von Neumann  
algebra acting on a separable 
Hilbert space $\HH$. Let $\A\subset\NN$ be a $*$-subalgebra, and 
${\mathcal M}\subset{\mathcal N}$ an $\A$-bimodule. Suppose that 
$\delta_1,...,\delta_k:\A\to{\mathcal N}$ are 
derivations such that the products 
$\delta_1(a_1)\cdots\delta_k(a_k)T\in {\mathcal M}$ 
for all $a_1,...,a_k\in\A$, where $T\in\NN$ is fixed. If 
$\phi:{\mathcal M}\to{\C}$ is a linear functional and we 
define $\tilde{\phi}\in C^k(\mathcal A)$ via:
\ben 
\tilde{\phi}(a_0,...,a_k)=\phi(a_0\delta_1(a_1)\cdots\delta_k(a_k)T),\een
then the Hochschild coboundary of $\tilde{\phi}$ is 
\ben 
(b\tilde{\phi})(a_0,...,a_{k+1})=(-1)^k\phi(a_0\delta_1(a_1)\cdots\delta_{k}(
a_k)a_{k+1}T-a_{k+1}a_0\delta_1(a_1)\cdots\delta_k(a_k)T).\een

\medskip

\proof{We prove this by induction. For $k=1$ we have
\bean (b\tilde{\phi})(a_0,a_1,a_2) &=& 
\phi(a_0a_1\delta(a_2)T-a_0\delta(a_1a_2)T+a_2a_0\delta(a_1)T)\nno
&=& -\phi(a_0\delta(a_1)a_2T-a_2a_0\delta(a_1)T).\eean
The derivation property shows that the first of these terms is still in 
${\mathcal M}$, so the $k=1$ case is true.
So we now suppose the result is true for all $n<k$. Then
\bean 
(b\tilde{\phi})(a_0,...,a_{k+1})&=&\phi(a_0a_1\delta_1(a_2)\cdots
\delta_k(a_{k
+1})T)\nno
&&\qquad +\sum_{i=1}^k(-1)^i\phi(a_0\delta_1(a_1)\cdots\delta_i(a_ia_{i+1})
\cdots
\delta_k(a_{k+1})T)\nno
&&\qquad +(-1)^{k+1}\phi(a_{k+1}a_0\delta_1(a_1)\cdots\delta_k(a_k)T)\nno
&=&(b\hat{\phi})(a_0,...,a_k)-(-1)^k\phi(a_ka_0\delta_1(a_1)\cdots\delta_{k-1}
(a_{k-1})\delta_k(a_{k+1})T)\nno
&&\qquad +(-1)^k\phi(a_0\delta_1(a_1)\cdots\delta_{k-1}(a_{k-1})\delta_k
(a_ka_{k+1})T\nno
&&\qquad +(-1)^{k+1}\phi(a_{k+1}a_0\delta_1(a_1)\cdots\delta_k(a_k)T).\eean
Here 
\ben 
\hat{\phi}(a_0,...,a_{k-1})=\phi(a_0\delta_1(a_1)\cdots\delta_{k-1}(a_{k-1})
(\delta_k(a_{k+1})T)).\een
Note that by hypothesis, the product of $\delta_1(a_1)\cdots
\delta_{k-1}(a_{k-1})$ and $\delta_k(a_{k+1})T$ is in ${\mathcal M}$. 

By induction we have
\bean 
(b\hat{\phi})(a_0,...,a_k)&=&(-1)^{k-1}\phi(a_0\delta_1(a_1)\cdots\delta_{k-1}
(a_{k-1})a_k\delta_k(a_{k+1})T)\nno
&&\qquad -(-1)^{k-1}\phi(a_ka_0\delta_1(a_1)\cdots\delta_{k-1}(a_{k-1})
\delta_k(a_{k+1})T).\eean
Thus we have
\ben 
(b\tilde{\phi})(a_0,...,a_{k+1})=(-1)^k\phi(a_0\delta_1(a_1)\cdots
\delta_k(a_k)a_{k+1}T-a_{k+1}a_0\delta_1(a_1)\cdots\delta_k(a_k)T).\een}
\end{lemma}

Thus the derivations need not all be the same to obtain the usual result 
linking Hochschild homology and derivations, \cite[p 84]{L}. We will mostly 
be interested in the case where the bimodule ${\mathcal M}$ is the 
ideal $\LL^{(1,\infty)}$, but we also use ${\mathcal M}=\LL^1$.

Our next aim is to define the Chern character of a finitely summable 
Fredholm module.
First we need a definition.

\begin{definition}
A {\bf pre-Fredholm module} for a unital Banach $*$-algebra
$\mathcal A$  is a pair $(\HH, F)$ where $\mathcal A$ is
(continuously) represented in $\mathcal N$ (a semifinite von Neumann
algebra acting on $\HH$) and $F$ is a self-adjoint
Breuer-Fredholm operator in $\mathcal N$ satisfying:

$1.\: 1-F^2 \in {\mathcal {K_N}},\:and$

$2.\: [F,a] \in {\mathcal {K_N}}\: for\: a \in {\mathcal A}.$

\noindent If $1-F^2=0$ we drop the prefix "pre-".

\noindent If, in addition, our module satisfies:

$1.'\: 1-F^2\in {\mathcal L}^{(p/2,\infty)}$

$2.'\: [F,a]\in {\mathcal L}^{(p,\infty)}$ 
for a dense set of $a\in \mathcal A$.
we say  $({\HH, F})$ is $(p,\infty)$-summable.
\end{definition}

{\bf Remark} Here and throughout the rest of the paper, if $p<2$ we interpret 
$T\in\LL^{(p/2,\infty)}$ as indicating that $T^{p/2}\in\LL^{(1,\infty)}$
which implies that $T\in\LL^1$.

Let  $(\HH,F)$ be a $p+1$-summable Fredholm module for $\A$ 
with $F^2=1$, 
\cite[IV.1.$\alpha$]{BRB}, that is, we have
$[F,a]\in\LL^{p+1}(\NN)$ for all $a\in\A$. In particular if
  $(\HH,F)$ is $(p,\infty)$-summable then it is $(p+1)$-summable.

The Chern character of $(\HH,F)$ is the class in periodic cyclic cohomology of 
the cocycles
\ben \lambda_n\tau'(\Gamma a_0[F,a_1]\cdots[F,a_n]), \ \ \ \ \ a_0,...,a_n
\in\A,\ \ \  n\geq p,\ \ \  n-p\ \ \mbox{even}.\een

Here $\lambda_n$ are constants ensuring that this collection of cocycles 
yields a well-defined periodic class, and they are given by
\ben \lambda_n=\left\{\begin{array}{ll} (-1)^{n(n-1)/2}\Gamma(\frac{n}{2}+1) 
& \ \ n\ \ \ {\rm even}\\ \sqrt{2i}(-1)^{n(n-1)/2}\Gamma(\frac{n}{2}+1) &
\ \ n\ \ \ {\rm odd}\end{array}\right..\een
The `conditional trace' (or, super-trace) $\tau'$ is defined by
\ben \tau'(T)=\frac{1}{2}\tau(F(FT+TF)),\een
provided $FT+TF\in\LL^1(\NN)$ (as it is in our case, see \cite[p293]{BRB}). 
Note that if $T\in\LL^1(\NN)$ 
we have (using the trace property and $F^2=1$)
\be \tau'(T)=\tau(T).\label{tracesequal}\ee
This class is represented by 
the cyclic cocycle $Ch_F\in C^p_\lambda(\A)$
\ben Ch_F(a_0,...,a_p)=\lambda_p\tau'(\Gamma a_0[F,a_1]
\cdots[F,a_p]), \ \ \ \ \ a_0,...,a_p\in\A.\een

If we only have  a pre-Fredholm module $(\HH,F)$, there is a canonical 
procedure described in \cite[p 310]{BRB} (and \cite{BF} in the general 
semifinite context) associating to $(\HH,F)$ a Fredholm module $(\HH',F')$. 
The Chern character of $(\HH,F)$ is then defined to be the Chern character 
of $(\HH',F')$. The Fredholm module $(\HH',F')$ has the same summability 
as $(\HH,F)$. We will not require the explicit form of this procedure, as 
we will now show that we have a more amenable procedure at our disposal.

Our next task is to show that if our spectral triple $(\A,\HH,\D)$
is such that $\D$ is not invertible, we can replace it by a new spectral 
triple in the 
same $K$-homology class in which the unbounded operator is invertible.
This is not a precise statement in the general semifinite case, as our 
spectral triples will not define $K$-homology classes in the usual sense. When 
we say that two spectral triples are in the same $K$-homology class, we shall 
take this to mean that the associated pre-Fredholm modules are operator 
homotopic up to the addition of degenerate Fredholm modules (see \cite{KK} 
for these notions, which make sense in our context).

\begin{definition} Let $(\A,\HH,\D)$ be a spectral triple. For any 
$m\in{\R}\setminus\{0\}$,  define the `double' of 
$(\A,\HH,\D)$ to be the spectral triple $(\A,\HH^2,\D_m)$ with 
$\HH^2=\HH\oplus\HH$, 
and the action of $\A$ and $\D_m$ given by
\ben \D_m=\bma \D & m\\ m & -\D\ema,\ \ \ \ a\to\bma a & 0\\ 0 & 0\ema,
\ \ \forall a\in\A.\een
\end{definition}

{\bf Remark} Whether $\D$ is invertible or not, $\D_m$ always is invertible, 
and $F_m=\D_m|\D_m|^{-1}$
has square 1. This is the chief reason for introducing this construction. We 
need to 
ensure that by doing so we do not alter the (co)homological data.

\begin{lemma}\label{noninv} The $K$-homology classes of $(\A,\HH,\D)$ and 
$(\A,\HH^2,\D_m)$ are the same. A representative of this class is 
$(\HH^2,F_m)$ 
with $F_m=\D_m|\D_m|^{-1}$.

\medskip

\proof{The $K$-homology class of $(\A,\HH,\D)$ is represented by the 
pre-Fredholm module 
$(\HH,F_\D)$ with $F_\D=D(1+\D^2)^{-1/2}$ while $[(\A,\HH^2,\D_m)]$ is 
represented 
by the 
pre-Fredholm module $(\HH^2,F_{\D_m})$ with $F_{\D_m}=\D_m(1+\D_m^2)^{-1/2}$ 
(we describe (pre)-Fredholm modules in subsection \ref{hochscyclic}). The one 
parameter 
family $(\HH,F_{\D_m})_{0\leq m\leq M}$ is a continuous operator homotopy, 
\cite{KK},\cite{CP1}, from 
$(\HH^2,F_{\D_M})$ to the 
direct sum of two pre-Fredholm modules

\ben (\HH,F_{\D})\oplus(\HH,-F_\D),\een
and in the odd case, the second pre-Fredholm module is operator homotopic to 
$(\HH,1)$ by the 
straight line path, since 
$\A$ is represented by zero on this module. In the even case 
we 
find the second pre-Fredholm module is homotopic to 
\ben \left(\HH,\bma 0 & 1\\ 1 & 0\ema\right),\een
the matrix decomposition being with respect to the ${\Z}_2$-grading of $\HH$.
Thus in both the even and odd cases the second module is degenerate, i.e. 
$F^2=1$, $F=F^*$ and $[F,a]=0$ for all $a\in\A$, and so 
the $K$-homology class of
$(\HH^2,F_{\D_M})$, written $[(\HH^2,F_{\D_M})]$,
is the K-homology class of 
$(\HH,F_{\D})$.
In addition, the Fredholm module $(\HH^2,F_m)$ with $F_m=\D_m|\D_m|^{-1}$ is 
operator 
homotopic to $(\HH^2,F_{\D_m})$ via
\ben t\to \D_m(t+\D_m^2)^{-1/2}\ \ \ 0\leq t\leq 1.\een
This provides the desired representative.}
\end{lemma}

The most basic consequence of Lemma \ref{noninv}, and the reason for proving 
it, comes from the following (see \cite[IV.1.$\gamma$]{BRB} and \cite{CIH}
for the proof).

\begin{proposition}\label{chind} The periodic cyclic cohomology class of the 
Chern character of a finitely summable Fredholm module depends only on its 
$K$-homology class.
\end{proposition}

In the general semifinite case this should be interpreted as saying that two 
pre-Fredholm modules which are operator homotopic up to the addition of 
degenerate Fredholm modules have the same Chern character.
In particular, therefore,
the Chern characters of $(\A,\HH,\D)$ and $(\A,\HH^2,\D_m)$ have the same 
class in periodic cyclic cohomology, and this can be computed using the 
Fredholm module $(\HH^2,F_m)$.

Using Connes' exact sequence, \cite[III.1.$\gamma$]{BRB},
\ben \cdots\sr{B}{\to}HC^p(\A)\sr{S}{\to}HC^{p+2}(\A)\sr{I}{\to}
HH^{p+2}(\A,\A^*)
\sr{B}{\to}HC^{p+1}(\A)\sr{S}{\to}\cdots\een
we see that the Hochschild class of $Ch_F$ is the image of the cyclic 
cohomology class of $Ch_F$ under the map $I$ induced 
by the inclusion of the cyclic complex in the Hochschild complex. This class 
is the 
noncommutative analogue of the integral representing the fundamental class. 
To see this,
recall that for the Dirac operator $\D$ on a closed spin manifold $X$ we have
\ben Ch(\D)(\cdot)=\mbox{const}\int_X \cdot\wedge\hat{A}=
\mbox{const}(\int_X\cdot\wedge 1
+ \int_X\cdot\wedge(-\frac{1}{24}p_1)+\cdots).\een
Here $\hat{A}$ is the A-roof genus, $p_i$ are the Pontryagin classes, and 
regarding 
$Ch(\D)$ as an element of de Rham homology, this formula tells us how to 
evaluate 
$Ch(\D)$ on elements of the exterior algebra of the manifold. In particular, 
restricting to differential forms of top degree (volume forms) we have
\ben Ch(\D)(f_0df_1\wedge\cdots\wedge df_{\dim X})=\mbox{const}\int_X 
f_0df_1\wedge\cdots\wedge df_{\dim X}.\een
Hence the Hochschild class of the Chern character yields the usual 
integration of a 
$(\dim X)$-form. This gives not only justification for the identification and 
study of
this Hochschild class, but also a heuristic for understanding the measurability
described in Corollary \ref{ismeas} (see Subsection \ref{31}).

Before leaving Chern characters, we note that the hypothesis of 
$(p,\infty)$-summability may be supplemented by Connes-Moscovici's 
discrete and finite dimension spectrum hypothesis, \cite{CM}. With this extra 
hypothesis 
one obtains a new representative of the Chern character expressed in terms of 
the 
operator $\D$. Using this representative, it is straightforward to identify 
the 
Hochschild class, and this agrees with the result stated 
in \cite[IV.2.$\gamma$]{BRB} 
and described here. However, the results concerning measurability (described 
later), arguably the most 
important consequence of Theorem \ref{thm8}, are rendered trivial, as the 
dimension 
spectrum hypothesis includes an assumption of measurability.

\subsection{The Dixmier Trace and the Heat Kernel}\label{dixheat}

Normally a Dixmier trace on the $\tau$-compact operators
means a
positive linear functional which is constructed
in the following way. One composes a positive
element $\omega$ of the dual of
$L^\infty({\R}^*_+)$ with the map which takes compact operators to the
Cesaro mean (described below) of their singular
values (where the latter is thought of as an element of
$L^\infty({\R}^*_+)$. The positive functional $\omega$ is also required to 
agree with the
ordinary limit on functions which have a limit at infinity.

The composition of any such $\omega$ from $L_\infty({\R}^*_+)^*$ which is 
vanishing on $C_0({\R}^*_+)$ 
with the Cesaro mean operator produces a functional which is (almost) 
dilation invariant and with which it is possible to define a 
non-normal trace (see [8]). We shall call such functionals Dixmier 
functionals and such non-normal traces Dixmier traces (see below).

A key technical lemma we will exploit uses the asymptotics of
the trace of the
heat operator for $\D$ to construct the singular or  Dixmier trace
that appears in Theorem \ref{thm8} when we have a particular kind of Dixmier 
functional $\omega$.

\begin{definition}
The Cesaro mean
on  $L^\infty({\R}^*_+)$, where ${\R}^*_+$ is the multiplicative
group of the positive reals, is given by:
$$M(g)(t)=\frac{1}{\log t}\int_1^t g(s) \frac{ds}{s}\;for\;g\in
L^\infty({\R}^*_+),\ t>0.$$
\end{definition}

\begin{definition}
We define the following maps on $L^\infty({\R}^*_+)$.
Let $D_a$ denote
dilation by $a\in {\R}^*_+$ and let $P^a$ denote exponentiation
by $a\in {\R}^*_+$. That is,
\begin{eqnarray*}
D_a(f)(x) &=& f(ax)\;for\;f\in L^\infty({\R}),and\\
P^a(f)(x) &=& f(x^a)\;for\;f\in L^\infty({\R}^*_+).
\end{eqnarray*}
$G$ is the set ${\R}^*_+ \times {\R}^*_+$ with multiplication:
$$ (s,t)(x,y) = (sx^t,ty).$$
\end{definition}

One of the main observations of \cite{BF} and \cite{CPS2} is that
in addition to dilation and Cesaro invariance, invariance under the 
operators $P^a$ ($a\in {\R}^*_+$)
is critical in one key step of the proof of the zeta function
representation of a Dixmier trace. We denote by $C_0({\R}^*_+)$ the 
continuous functions
on ${\R}^*_+$ vanishing at infinity. We will need the 
existence of a $G$-invariant, $M$ invariant Dixmier functional
 on $L^\infty({\R}^*_+)$.

\begin{theorem}[\cite{CPS2}]\label{niceomega}
There exists a state $\Omega$ on $L^\infty({\R}^*_+)$
satisfying the following conditions:\\
(1) $\Omega(C_0({\R}^*_+)) \equiv 0$.\\
(2) If $f$ is real-valued in $L^\infty({\R}^*_+)$
then
$$ess\ lim\mbox{-}inf_{t\to\infty} f(t) \leq\Omega(f)\leq ess\ lim\mbox{-}sup_{t\to\infty}
f(t).$$
(3) If the essential support of $f$ is compact then $\Omega(f)=0.$\\
(4) For all $c\in {\R}^*_+$, $\Omega(D_cf)=\Omega(f)$ for all
$f\in L^\infty({\R}^*_+)$.\\
(5) For all $a \in{\R}^*_+$ and all $f\in L^\infty({\R}^*_+)$
$\Omega(P^af)=\Omega(f)$.\\
(6) For all  $f\in L^\infty({\R}^*_+)$, $\Omega(Mf)=  \Omega(f)$.
\end{theorem}

The approach of \cite{CPS2} as described in Theorem \ref{niceomega} is to 
construct
what might be more appropriately be termed a `maximally
invariant Dixmier functional'. This maximal invariance
is what is required to establish the zeta function
representation of a Dixmier trace (and hence the heat kernel
formula for $\LL^{(1,\infty)}$) in full generality. Weaker conditions suffice
for  the case of $L^{(p,\infty)}$, $p>1$,  essentially because
the map $T\to T^p$ taking $L^{(p,\infty)}$ to $L^{(1,\infty)}$
is not surjective and in fact the image is a smaller ideal consisting of
compact operators $T$ whose singular values satisfy, for some $C>0$,
the inequality
$\mu_s(T) \leq C/s$
for $s$ suficiently large; see \cite{CPS2} for further discussion.

A notation we will often use is to write, for a given function
$f\in L^\infty({\R}^*_+)$ and Dixmier functional $\omega$,
$\omega(f) =\omega\mbox{-}\hspace{-.03in}\lim_{\lambda\to\infty}f(\lambda)$. In particular we 
will be interested in applying such functionals to the function
\ben \frac{1}{\log(1+t)}\int_0^t\mu_s(T)ds\een
where $T\in \LL^{(1,\infty)}$ is positive. This is the Dixmier trace 
associated to the semifinite normal trace $\tau$, denoted $\tau_\omega$, and 
we extend it to all of $\LL^{(1,\infty)}$ by linearity. The Dixmier 
trace $\tau_\omega$ is defined on the ideal $\LL^{(1,\infty)}$, and vanishes 
on the ideal of trace class operators. This latter fact is used repeatedly 
throughout the paper without further comment.

Let   $T\geq 0$ and 
define $e^{-T^{-2}}$ as the operator that is zero on $\ker T$
and on $\ker T^\perp$ is defined in the usual way by the functional
calculus. We remark that if $T\geq 0$, $T\in{\mathcal L}^{(p,\infty)}$
for some $p\geq 1$ then  $e^{-tT^{-2}}$ is trace class for all $t>0$.
Then we have

\begin{theorem}[\cite{CPS2}]\label{cps2}
If $A\in\mathcal N$, $T\geq 0$, $T\in{\mathcal L}^{(p,\infty)}$
then,
$$\Omega\mbox{-}\hspace{-.05in}\lim_{\lambda\to\infty}\lambda^{-1}\tau(Ae^{-\lambda^{-2/p}T^{-2}})=
\Gamma(p/2+1)\tau_\Omega(AT^{p})$$
for $\Omega\in L^\infty({\R}^*_+)^*$ satisfying the conditions of Theorem 
\ref{niceomega}.
\end{theorem}

{\bf Remark}. The reason for the citation of \cite{CPS2} for this result is 
that we
require the case $p=1$, and this is the only place where this is established. 
For $p>1$, however, see \cite[p563]{BRB} and \cite{Va}.

To use this result in this paper we will apply it to the case
where $T=(1+\D^2)^{-1/2}$ or $T=\dd^{-1}$ if $\D$ has bounded inverse. Then
a simple but useful corollary of this theorem is that for $p\geq 1$ 
and $\dd^{-1}\in{\mathcal L}^{(p,\infty)}$
the function on ${\R}^*_+$ given by 
$$t \to t^p\tau(Ae^{-t^2\D^2})$$ is bounded. This follows from setting
$\lambda^{-1}=t^p$ and $T=\dd^{-1}$. Or in other words

\begin{lemma}\label{heat} If $p\geq 1$ and $(\A,\HH,\D)$ is a 
$(p,\infty)$-summable spectral triple with $\D$ invertible, then there exists 
a constant $C_p>0$ 
such that
\ben \tau(e^{-t^2\D^2})\leq C_pt^{-p}\ \ \mbox{for}\ \ t> 0.\een
\end{lemma}

A fact that we will frequently require is the following.

\begin{proposition}[\cite{CPS2,It,FN}]\label{hypertrace} 
The Dixmier trace $\tau_\omega$ associated to a Dixmier functional $\omega$\\ 
defines a trace on
the algebra of a 
$QC^1$ $(p,\infty)$-summable spectral triple via
\ben a\mapsto \tau_\omega(a(1+\D^2)^{-p/2}).\een
\end{proposition}

\section{The Hochschild Class of the Chern Character}\label{MR}

\subsection{Statement of the Main Result}\label{31}
Our main result is the general semifinite version of 
a Type I result in \cite[IV.2.$\gamma$]{BRB} 
which identifies the Hochschild class of 
the Chern 
character of a $(p,\infty)$-summable spectral triple. 
With the preliminary definitions out of the way, we can now state our main 
result:

\begin{theorem}\label{thm8} Let $(\A,\HH,\D)$ be a $QC^k$ 
$(p,\infty)$-summable spectral triple with $p\geq 1$ integral and 
$k=\max\{2,p-2\}$. Then

1) A Hochschild cocycle on $\A$ is defined by
\ben \phi_\omega(a_0,...,a_p)=\lambda_p\tau_\omega(\Gamma 
a_0[\D,a_1]\cdots[\D,a_p](1+\D^2)^{-p/2}),\een

2) For all Hochschild $p$-cycles $c\in C_p(\A)$ (i.e., $bc=0$), 
\ben\la \phi_\omega,c\ra=\la Ch_{F_\D},c\ra,\een
where $Ch_{F_\D}$ is the Chern character in cyclic cohomology of the 
pre-Fredholm module 
over $\A$ with $F_\D=\D(1+\D^2)^{-1/2}$. \end{theorem}

{\bf Remark} Here $\tau_\omega$ is the Dixmier trace associated to any Dixmier 
functional $\omega$.
The two most important corollaries of Theorem \ref{thm8} are the following.
\begin{corollary}\label{ismeas} Let $(\A,\HH,\D)$ be as in Theorem \ref{thm8}. 
If $c=\sum_i a_0^i\otimes a_1^i\otimes\cdots \otimes 
a_p^i$ is a Hochschild $p$-cycle, then 
\ben \Gamma\sum_ia_0^i[\D,a_1^i]\cdots[\D,a_p^i](1+\D^2)^{-p/2}\een
is measurable.

\end{corollary}

{\bf Remark} An operator $T\in\LL^{(1,\infty)}$ is measurable (in the sense of
Connes) if the 
$\omega$-limit
$$\omega\mbox{-}\hspace{-.05in}\lim_{t\to\infty}\frac{1}{\log(1+t)}\int_0^t\mu_s(T)ds$$
is independent of the choice of $\omega$. We will include a proof of this 
important result (Corollary 11) as part of the proof of Theorem \ref{thm8}.

\begin{corollary} With $(\A,\HH,\D)$ as in Theorem \ref{thm8}, and supposing 
that $Ch_{F_\D}$ pairs nontrivially with $HH_p(\A)$, then
\ben \tau_\omega((1+\D^2)^{-p/2})\neq 0.\een
\end{corollary}

{\bf Remark} The hypothesis of the Corollary is that there exists some 
Hochschild $p$-cycle such that $\la ICh_{F_\D},c\ra\neq 0$. Computing this 
pairing using Theorem \ref{thm8} above, we see that  $(1+\D^2)^{-p/2}$ can 
not have zero Dixmier trace for any  choice of Dixmier functional $\omega$. 
For if $(1+\D^2)^{-p/2}$ did have vanishing Dixmier trace, 
and $c=\sum_ia^i_0\otimes\cdots\otimes a^i_p$ is any Hochschild cycle
\bean |\la ICh_{F_\D},c\ra|&=&\left|\sum_i\tau_\omega
\left(\Gamma a^i_0[\D,a^i_1]\cdots[\D,a^i_p](1+\D^2)^{-p/2}\right)
\right|\nno   
&\leq&\sum_i\n\Gamma a^i_0[\D,a^i_1]\cdots[\D,a^i_p]\n\tau_\omega
\left((1+\D^2)^{-p/2}\right)=0.\eean
Hence if the pairing is nontrivial, the Dixmier trace can not vanish 
on $(1+\D^2)^{-p/2}$.

During the course of the proof {\em we will always suppose} that we have a 
spectral 
triple $(\A,\HH,\D)$ with $\D$ invertible, by replacing $(\A,\HH,\D)$ by 
$(\A,\HH^2,\D_m)$ if necessary. Despite knowing that the cyclic classes of 
the 
Chern characters of these two triples coincide, by Lemma \ref{noninv} and 
Proposition \ref{chind}, and so their Hochschild classes also coincide, we do 
not know that this is true for 
the specific representative displayed in Theorem \ref{thm8}, and this is 
something we will need to determine. 
A proof that this is indeed the case can be found in the Appendix.

Before discussing the proof any further, we show that the functional 
$\phi_\omega$ is indeed a Hochschild cocycle.

\begin{lemma}\label{iscocycle} Let $p\geq 1$ and suppose that $(\A,\HH,\D)$ 
is a $QC^1$ 
$(p,\infty)$-summable spectral triple. Then the multilinear functional
\ben \phi_\omega(a_0,...,a_p)=
\lambda_p\tau_\omega(\Gamma a_0[\D,a_1]\cdots[\D,a_p](1+\D^2)^{-p/2})\een
is a Hochschild cocycle.

\medskip

\proof{By Lemma \ref{hochs} and the trace property of the Dixmier trace, we 
have
\bean (b\phi_\omega)(a_0,...,a_p)&=&(-1)^{p-1}\lambda_p
\tau_\omega(\Gamma a_0[\D,a_1]
\cdots[\D,a_{p-1}]a_p(1+\D^2)^{-p/2})\nno
&&\qquad-(-1)^{p-1}\lambda_p\tau_\omega(\Gamma a_0[\D,a_1]
\cdots[\D,a_{p-1}](1+\D^2)^{-p/2}a_p).\eean
As $(\A,\HH,\D)$ is $QC^1$, 
\bean [(1+\D^2)^{-p/2},a_p]
&=&-\sum_{k=0}^{p-1}(1+\D^2)^{-(p-k)/2}[(1+\D^2)^{1/2},a_p]
(1+\D^2)^{-(1+k)/2},\eean
and this is trace class. So  $a_p(1+\D^2)^{-p/2}=(1+\D^2)^{-p/2}a_p$ 
modulo trace class operators, and so
the two terms above cancel.}
\end{lemma}

Thus to show that $\phi_\omega$ is a Hochschild cocycle is relatively simple, 
and does 
not require the full smoothness assumptions of Theorem \ref{thm8}. Of course 
the 
important aspects of Theorem \ref{thm8} are that $\phi_\omega$ is a 
representative of 
the Hochschild class of the Chern character, and the measurability 
of `$p$-forms'.

\subsection{What was previously known}\label{32}
This theorem, for $\NN=\B(\HH)$ and $1<p<\infty$ ($p$ integral) was proved in 
lectures by Alain Connes at the Coll\`{e}ge de France in 1990. A version of 
this argument appeared in \cite{Va}. The extension of this argument to 
general semi-finite von Neumann algebras, with the additional hypothesis 
that $\D$ have bounded inverse, is presented in the preprint of 
Benameur-Fack, \cite{BF} and we thank the authors for bringing it to our
attention. It provided an impetus to our work.
Some supplementary details in the proof were 
given to us by Thierry Fack, and we thank him for his notes, \cite{FN}. 
In addition, a simpler strategy
using the pseudodifferential calculus of Connes-Moscovici, 
\cite{CM}, was communicated to us by Nigel Higson. 
In conjunction with the results in 
\cite{CPS2}, Higson's argument appears to generalise to the semifinite case as 
well as giving an alternate proof of Theorem \ref{thm8}, however 
we will not describe the details
here.

The 
extension of these earlier results which our Theorem \ref{thm8} implies are

1) for the first time we provide a proof for the case $p=1$ (the proof in 
this case overcomes some serious technical obstacles).

2) We dispense with the hypothesis in the type $II_\infty$ case that $\D$ 
has bounded inverse. This is crucial due to the
`zero-in-the-spectrum'  
phenomenon for $\D$ (that is, for type II $\NN$, zero is generically in the 
point and/or continuous spectrum, \cite{FW}) and is not just the simple
problem posed by non-trivial $\ker \D$.

3) Importantly, our strategy of proof is the same for all $p\geq 1$, 
is independent of the type of the von Neumann algebra $\NN$ and is simpler
than previously published arguments.

We now come to the proof of Theorem \ref{thm8}. 
The general form of the technical estimates, and so the basic structure of the 
analytic parts of the proof, are based on a synthesis of our understanding of 
the 
arguments in \cite{BF} and \cite{Va}. These in turn have their origin in the 
original 
arguments of Connes. The latter parts of the argument where we need to 
construct various
cohomologies in the Hochschild theory to arrive at the functional in the 
statement of 
the theorem, closely follow the argument in \cite{Va}.

Our method of proof is in some ways more direct than these other approaches. 
In particular we do not 
need to prove our technical estimates for general functions of $\D$, only 
the particular functions that allow us to employ the heat kernel approach to 
the Dixmier trace. The chief novelty (and difficulty) of this direct 
approach is that we can deal with the case $p=1$. For this approach the 
assumption of \cite{Va} that the functions 
of $\D$ involved are 
compactly supported is of no use and
various technical estimates in \cite{BF,BRB,Va} are not available.

\subsection{Functional Calculus Preliminaries}

In this subsection we establish some trace and commutator estimates for 
certain 
functions of $\dd$. We will work exclusively with one function, however the
 definition 
of this function depends on the value of $p$. Moreover there are substantial 
differences 
between the even and odd cases, and for technical reasons we also require 
estimates 
involving square roots of functions.

For $p\geq 1$ an integer, and $x\geq 0$
define
\be erf_p(x)=\frac{p}{\Gamma(\frac{p}{2}+1)}\int_0^x 
r^{p-1}e^{-r^2}dr.\label{erf}\ee
Using
\ben \int_0^\infty r^{p-1}e^{-r^2}dr=\frac{\Gamma(p/2)}{2}=\frac{\Gamma
(\frac{p}{2}+1)}{p},\een
we have $erf_p(\infty)=1$ and $erf_p(0)=0$. Now define
\be f_p(x) =\left\{\begin{array}{ll} 1-erf_p(x) & x\geq 0\\
1-(-1)^perf_p(-x) & x\leq 0\end{array}\right.\label{fprov}\ee

Then we have $f_p(0)=1$, $f_p(\infty)=0$ and 
\be 
f'_p(x)=\frac{-p}{\Gamma(\frac{p}{2}+1)}x^{p-1}e^{-x^2}.\label{fprime}\ee
For $p$ even and all $x\in{\R}$ or $p$ odd and $x\geq 0$ we can write 
\bea f_p(x)&=& 1-erf_p(|x|)\nno 
&=& \frac{p}{\Gamma(p/2+1)}\int_1^\infty x^ps^{p-1}e^{-s^2x^2}ds\nno
&=& c(p)\int_1^\infty x^ps^{p-1}e^{-s^2x^2}ds.\label{better}\eea

We will see shortly that $f_p$ is Schwartz class for $p$ even. However for 
$p$ odd, 
while $f_p(x)\to 0$ rapidly as $x\to +\infty$, as $x\to-\infty$, 
$f_p(x)\to 2$. The 
reason we have defined the function in this way is to obtain smoothness at 
$x=0$, and 
the important part of the definition is for $x\geq 0$ anyway. For instance, we 
have our 
first estimate.

\begin{lemma}\label{intable} Let $p\geq 1$ and suppose that $(\A,\HH,\D)$ is 
a $QC^0$ $(p,\infty)$-summable spectral triple with $\D$ invertible. If $h$ is 
either of the functions $f_p$ or $\sqrt{f_p}$ then for $t>0$
\ben \n h(t\dd)\n_1\leq C_ht^{-p}.\een

\medskip

\proof{Let $d\phi_\lambda=d\tau(E_\lambda)$ be the scalar spectral measure 
for $\dd$, and consider 
the function $\sqrt{f_p}$. We have by Lemma 8.2 of \cite{CP2}
\bean \tau(\sqrt{f_p}(t\dd)) & =& 
\left(\frac{p}{\Gamma(\frac{p}{2}+1)}\right)^{1/2}\tau\left(\left(
\int_1^\infty s^{p-1}t^p\dd^pe^{-s^2t^2\D^2}ds\right)^{1/2}\right)\nno
&=&\left(\frac{t^pp}{\Gamma(\frac{p}{2}+1)}\right)^{1/2}\int_0^\infty
\left(\int_1^\infty 
s^{p-1}\lambda^pe^{-s^2t^2\lambda^2}ds\right)^{1/2}d\phi_\lambda\nno
&\leq & \left(\frac{t^pp}{\Gamma(\frac{p}{2}+1)}\right)^{1/2}\int_0^\infty 
\lambda^{p/2}e^{-t^2\lambda^2/4}\left(\int_1^\infty 
s^{p-1}e^{-s^2t^2\lambda^2/2}ds\right)^{1/2}d\phi_\lambda\nno
&\leq & \left(\frac{t^pp}{\Gamma(\frac{p}{2}+1)}\right)^{1/2}\int_0^\infty 
\lambda^{p/2}e^{-t^2\lambda^2/4}\left(\int_0^\infty 
s^{p-1}e^{-s^2t^2\lambda^2/2}ds\right)^{1/2}d\phi_\lambda\nno
&=&\left(\frac{t^pp}{\Gamma(\frac{p}{2}+1)}\right)^{1/2}\int_0^\infty 
\lambda^{p/2}e^{-t^2\lambda^2/4}(\Gamma(p/2)
\lambda^{-p}t^{-p}2^{\frac{p}{2}-1})^{1/2}
d\phi_\lambda\nno
&=& 2^{\frac{p}{4}}\int_0^\infty e^{-t^2\lambda^2/4}
d\phi_\lambda\nno
&=& 2^{\frac{p}{4}}\tau(e^{-(t/2)^2\D^2})\nno
&\leq & Ct^{-p},\eean
where the last line follows from the heat kernel estimate Lemma \ref{heat}. 
The same method applies to yield the result for $h=f_p$ also.}
\end{lemma}

The above Lemma required knowledge of $f_p$ for positive arguments, so for $p$ 
odd, 
we are free to alter the definition in any reasonable way for negative values.

So for $p\geq 1$ odd, and for some $k>0$, define
\ben f_p(x)=\left\{\begin{array}{ll} 1-erf_p(x) & x\geq 0\\
1+erf_p(-x) & -k\leq x\leq 0\\
g(x) & x\leq -k\end{array}\right.\een
Here we take $g(x)=Q(x)e^{-x^2}$, where $Q$ is a polynomial. We may choose to 
make 
$f_p$ a $C^l$ function at $-k$, and this will require taking $Q$ to be of 
order $l$.

\begin{lemma} For $p$ even
\ben f_p(x)=\sum_{i=0}^{[p-2/2]}\frac{c(p-2i)}{2}x^{p-2-2i}e^{-x^2}\een
When $p$ is odd and $p\geq 3$, and $x\geq 0$ we have 
\ben f_p(x)=\sum_{i=0}^{[p-2/2]}
\frac{c(p-2i)}{2}x^{p-2-2i}e^{-x^2}+f_1(x).\een

\medskip
\proof{Integration by parts using the formulae in Equations (\ref{better}), 
and the observation that for all integers $p\geq 3$
\ben \frac{c(p)}{c(p-2)}\frac{(p-2)}{2}=1.\een}
\end{lemma}

Observe that the Lemma shows that for $p$ even, the function $f_p$ is Schwartz 
class. 
Indeed, $f_p^{1/2}$ is Schwartz class. This follows because 
$f_p(x)=P(x)e^{-x^2}$ where 
$P$ is an even polynomial with a nonzero constant term. Thus
\be \frac{df_p^{1/2}}{dx}=\frac{(P'(x)-2xP(x))e^{-x^2/2}}{P(x)^{1/2}},
\label{dies}\ee
from which it is easy to see that $f_p^{1/2}$ has derivatives of all orders, 
and they 
are all of rapid decrease.

For  $p$ odd we will have a rapidly decaying function also, but only $C^l$ at 
$-k$, 
where we may choose $l$ as large as we like. To see this for large positive 
$x$ we 
require the following result, \cite{Dwi}.

\begin{lemma}\label{rapid} There is an asymptotic expansion of 
\ben f_1(x)=\frac{2}{\sqrt{\pi}}\int_{x}^\infty e^{-s^2}ds\een
as $x\to +\infty$ of the form
\ben f_1(x)\sim \frac{e^{-x^2}}{\sqrt{\pi}x}\left[1-\frac{2!}{(2x)^2}+
\frac{4!}{2!(2x)^4}-\frac{6!}{3!(2x)^6}+\cdots\right].\een

\medskip

\proof{The function $f_1$ is precisely the complementary error function 
(for positive $x$), and there is a standard asymptotic expansion for $x$ large 
and
positive
\ben erfc(x)=f_1(x)\sim \frac{e^{-x^2}}{\sqrt{\pi}x}
\left[1-\frac{2!}{(2x)^2}+\frac{4!}
{2!(2x)^4}-\frac{6!}{3!(2x)^6}+\cdots\right].\een}
\end{lemma}

At this point we know enough to proceed when $p$ is even, but for the 
estimates we wish 
to prove next, we require more information for $p$ odd.

For $p$ odd, our definition of $f_p$ ensures that $x^mf_p^{1/2}$ is 
integrable for all 
$m\geq 0$, so the Fourier transform of $f_p^{1/2}$ is smooth, and of course 
lies in 
$C_0({\R})$. If we define $f_p$ so that it is  $C^l$, then the first $l$ 
derivatives of 
$f_p^{1/2}$ will also have smooth Fourier transform, contained in 
$C_0({\R})$, using 
Lemma 15 and an argument similar to that in Equation \ref{dies}. So 
for 
$l\geq i+2$, the Fourier transform of $\p^if_p^{1/2}$ is in $L^1({\R})$. This 
follows 
because
\ben \widehat{(\p^{i+2}f_p^{1/2})}(\xi)\to 0\ \ \mbox{as}
\ \ |\xi|\to\infty,\een
so
\ben |\xi^2\widehat{(\p^if_p^{1/2})}(\xi)|\to 0\ \ \mbox{as}
\ \ |\xi|\to\infty,\een
which tells us that
\ben |\widehat{(\p^if_p^{1/2})}(\xi)|=o(|\xi|^{-2}) \ \ \mbox{as}
\ \ |\xi|\to \infty.\een
Choosing $l\geq 4$ then tells us that in both the even and odd cases
\be \int_{\R}|\xi^i\widehat{f_p^{1/2}}(\xi)|
d\xi<\infty,\ \ i=1,2.\label{i12}\ee

We use this to formulate two commutator estimates.

\begin{lemma}\label{1bound} If $(\A,\HH,\D)$ is a 
$QC^1$ $(p,\infty)$-summable spectral 
triple then
\ben \n [f_p(t\dd),a]\n_1\leq C_{a,f_p,p}t^{-p+1}.\een

\medskip

\proof{Writing $h=f_p^{1/2}$, we have the straightforward calculation
\bean \n [f_p(t\dd),a]\n_1 &=&\n h(t\dd)[h(t\dd),a]+[h(t\dd),a]h(t\dd)\n_1\nno
&\leq & 2\n[h(t\dd),a]\n_\infty\n h(t\dd)\n_1\nno
&\leq & 2C_pt\n[\dd,a]\n_\infty\n\widehat{h}(\xi)\xi\n_1t^{-p}\nno
&=& C_{a,f_p,p}t^{-p+1}.\eean
The last inequality comes from Lemma \ref{intable}
and 
\bean [h(t\dd),a] &=& \frac{1}{\sqrt{2\pi}}\int_{\R}\hat{h}(s)[e^{its\dd},a]
ds\nno
&=&\frac{1}{\sqrt{2\pi}}\int_{\R}\hat{h}(s)is\int_0^1e^{itsr\dd}[t\dd,a]
e^{i(1-r)st\dd}drds.\eean
The finiteness of (\ref{i12}) for $i=1$ completes the proof.}
\end{lemma}

\begin{lemma}\label{tracechainrule} If $(\A,\HH,\D)$ is a 
$QC^2$ $(p,\infty)$-summable 
spectral triple then
\ben \n [f_p(t\dd),a]-\frac{1}{2}\{f'_p(t\dd),t[\dd,a]\}\n_1\leq 
C_{a,p,f_p}t^{-p+2},\een
where $\{T,S\}=TS+ST$.

\medskip
\proof{Again we write $h=f_p^{1/2}$, and again this is just a computation.
\bean & &\n [f_p(t\dd),a]-\frac{1}{2}\{f_p'(t\dd),t[\dd,a]\}
\n_1\qquad\qquad\qquad\nno
&=& \n h(t\dd)[h(t\dd),a]+[h(t\dd),a]h(t\dd)\nno 
&&\qquad-  h(t\dd)h'(t\dd)t[\dd,a]-t[\dd,a]h(t\dd)h'(t\dd)\n_1\nno
&\leq & \n [h(t\dd),a]-th'(t\dd)[\dd,a]\n_\infty\n h(t\dd)\n_1\nno 
&&\qquad+ \n[h(t\dd),a]-t[\dd,a]h'(t\dd)\n_\infty\n h(t\dd)\n_1\nno
&\leq & t^{-p+2}C_p\int_{\R}|\widehat{h}(\xi)\xi^2|d\xi.\eean
The final inequality follows from Lemma \ref{intable} and writing, 
\cite{BF} 
\ben A(t):=[h(t\dd),a]-h'(t\dd)[t\dd,a],\een
we have
\bean A(t) &= 
&\int_{\R}\hat{h}(u)\int_0^1(e^{iuts\dd}[iut\dd,a]e^{iut(1-s)\dd}-e^{iut\dd}
[iut\dd,a])dsdu\nno
&=&\int_{\R}\hat{h}(u)\int_0^1
e^{iuts\dd}[[iut\dd,a],e^{iut(1-s)\dd}]dsdu\nno
&=& 
-\int_{\R}\hat{h}(u)\int_0^1e^{iuts\dd}\int_0^1e^{iut(1-s)r\dd}[iut(1-s)\dd,
[iut\dd,a]]e^{iut(1-s)(1-r)\dd}dsdrdu.\eean
A similar result holds for $B(t)=[h(t\dd),a]-[t\dd,a]h'(t\dd)$. In both cases 
$\n A(t)\n_\infty$ and $\n B(t)\n_\infty$ are $O(t^2)$ as $t\to 0$. The 
finiteness of 
(\ref{i12}) for $i=2$ completes the proof.}
\end{lemma}

Estimates like those presented in the last two Lemmas may be regarded as 
approximate 
extensions of familiar rules of calculus to the `quantum' setting. Both the 
previous 
Lemmas extend to a large class of functions, but as we only require these very 
particular results, we do not pursue these matters here.

\subsection{From the Chern Character to a Hochschild Cocycle}

Now that we know something about $f_p$, we can begin the proof. The first step 
is to 
bring $\dd$ into the picture, and to do so in a way that will allow us, 
eventually, to 
make use of its summability.

\begin{lemma}[\cite{BF,Va}]\label{lim} Let $p\geq 1$ be integral and suppose 
that $(\A,\HH,\D)$ is 
a $QC^0$ $(p,\infty)$-summable spectral triple with $\D$ invertible.
Let $c=\sum_ia^i_0\otimes\cdots\otimes a^i_p$ 
be a Hochschild $p$-cycle.  Then
\ben \la ICh^*(F),c\ra=-\lim_{t\to 0}\lambda_p\sum_i\tau(\Gamma 
a^i_0[F,a^i_1]\cdots[F,a^i_{p-1}]F[f_p(t\dd),a^i_{p}]),\een
where  $F_\D=D\dd^{-1}$.

\medskip

\proof{Ignore $i$ momentarily, and set $A=\Gamma a_0[F,a_1]\cdots[F,a_p]$. 
As $t\to 0$ we have $f_p(t\dd)\to 1$ (strong operator topology), so  
\ben \la ICh^*(F_\D),c\ra=\lim_{t\to 0}\tau'(f_p(t\dd)A).\een
Here we have used $Ff_p(t\dd)=f_p(t\dd)F$ to see that the right hand side is 
equal to
\ben \tau(f_p(t\dd)F(FA+AF)).\een
As $F(FA+AF)$ is trace class, \cite[I.6.1,p93]{Dix2} shows that the above 
equality holds.

For $t>0$, the operator $f_p(t\dd)$ is trace class by Lemma \ref{intable}, 
so we may replace $\tau'$ by $\tau$, using Equation \ref{tracesequal}. Making 
this
change and expanding the 
last factor of $A$ gives
\ben \tau'(f_p(t\dd)A)
=\tau(\Gamma a_0[F,a_1]\cdots[F,a_{p-1}]Fa_pf_p(t\dd))-\tau(\Gamma 
a_0[F,a_1]\cdots[F,a_{p-1}]a_pFf_p(t\dd)).\een

Using the fact that $a\to [F,a]$ is a derivation, we can use Lemma 
\ref{hochs} and
\ben \sum_ib(a_0^i\otimes a_1^i\otimes\cdots\otimes a_p^i)=0,\een
to see that
\ben \tau(\Gamma a_0[F,a_1]\cdots[F,a_{p-1}]a_pFf_p(t\dd))-\tau(\Gamma 
a_pa_0[F,a_1]\cdots[F,a_{p-1}]Ff_p(t\dd))=0,\een
as this is a Hochschild coboundary paired with a Hochschild cycle. 
This proves the Lemma.}
\end{lemma}

Note this {\em only} works when we pair with a Hochschild cycle. For an 
arbitrary 
chain we can not swap $a_p$ around to the front. 
Nevertheless, for any $a_0,...,a_p\in\A$ we can define a one-parameter family 
of 
multilinear functionals
\ben \psi_t(a_0,...,a_p):= -\lambda_p\tau(\Gamma 
a_0[F,a_1]\cdots[F,a_{p-1}]F[f_p(t\dd),a_p]),\een
and we have already shown that the pairing of the 
Chern character with Hochschild cycles is given by pairing the Hochschild 
cycle with 
$\psi_t$ and taking the limit as $t\to 0$. 
However, we {\em have not} yet defined a multilinear functional which 
represents $ICh_F$. {\em If} we knew that for all $a_0,...,a_p\in\A$
\ben \phi(a_0,...,a_p):=\lim_{t\to 0}\psi_t(a_0,...,a_p)\een
existed, and we could show that $b\phi=0$, then we would have
\ben [\phi]=[ICh_F]\in HH^*(\A,\A^*).\een
In general we can not assert the existence of the above limit, and this is 
why we do not yet have a representative of the Hochschild class of the Chern 
character.

The strategy is to show that $|\psi_t(a_0,...,a_p)|$ is {\em bounded} as 
$t\to 0$ for any $a_0,...,a_p\in\A$, so that we may define a functional by 
taking the $\omega$-limit. We will then rewrite this result in terms of the 
associated Dixmier trace. Once achieved, we will have a well-defined 
multilinear functional on $\A$ which depends on the choice of $\omega$. 
However, the pairing of this functional with Hochschild cycles will return 
the true limit, no matter what choice of $\omega$ is employed. This is the 
origin of Corollary \ref{ismeas}, but the precise form must wait until we 
have identified $\phi_\omega$ as a representative of the Hochschild class.

So to begin, let us obtain the estimate which will allow us to show 
that $\psi_t$ 
is bounded.

\begin{lemma}\label{p1bound} Let $p\geq 1$ and suppose that $(\A,\HH,\D)$ is 
a $QC^1$ $(p,\infty)$-summable spectral triple with $\D$ invertible. Then 
\ben \n[f_p(t\dd),a]\n_{(p,1)}\ \ {\rm is\ bounded\ as}\ t\to 0.\een

\medskip

\proof{We have the estimate, Lemma 
\ref{1bound}, for all $p\geq 1$
\ben \n[f_p(t\dd),a]\n_1 \leq  \tilde{C}_ft^{-p+1}\n[\dd,a]\n_\infty.\een
So for $p=1$ we 
are done. For $p>1$ we have the interpolation inequality
\ben \n T\n_{(p,1)}\leq C_p\n T\n_1^{1/p}\n T\n_\infty^{1-1/p},\ \ 
T\in\LL^1.\een

In Lemma \ref{1bound} we also estimated the norm, obtaining
\ben \n[f_p(t\dd),a]\n_\infty=O(t)\een
which allows us to finish the proof since
\ben \n[f_p(t\dd),a]\n_{(p,1)}\leq {\cal 
C}_{f_p,p}(t^{-p+1})^{1/p}t^{1-1/p}={\cal C}_{f_p,p}.\een}
\end{lemma}

{\bf Remark} The use of the interpolation inequality in the previous proof is 
standard in the type I setting, \cite[IV, Appendix B]{BRB}. For the type II 
case we note that it is sufficient to obtain the result for the commutative 
von Neumann algebra $L^\infty(0,\infty)$ and apply the results of \cite{DDP}.

\begin{lemma}\label{psitozero} Let $p\geq 1$ and suppose that $(\A,\HH,\D)$  
is a $QC^1$ 
$(p,\infty)$-summable spectral triple with $\D$ invertible. 
Then for all $a_0,...,a_p\in\A$ the 
function
\ben t\to \psi_t(a_0,...,a_p)\een
is bounded as $t\to 0$.

\medskip
\proof{By Lemma \ref{sizes}, $[F,a_i]\in\LL^{(p,\infty)}$, $i=1,...,p$. So 
\ben \Gamma a_0[F,a_1]\cdots[F,a_{p-1}]F\in\LL^{(q,\infty)},\een
where $q=p/(p-1)$ (for $p=1$ replace the $(q,\infty)$ norm with the operator 
norm). The K\"{o}the dual of $\LL^{(p,1)}$ is 
$\LL^{(q,\infty)}$, so as $t\to 0$
\bean |\psi_t(a_0,...,a_p)|&\leq &\n\Gamma 
a_0[F,a_1]\cdots[F,a_{p-1}]F\n_{(q,\infty)}\n[f_p(t\dd),a_p]\n_{(p,1)}\nno
&\leq & {\cal C}_{f_p,p}\n\Gamma 
a_0[F,a_1]\cdots[F,a_{p-1}]F\n_{(q,\infty)},\eean
by Lemma \ref{p1bound}.}
\end{lemma}

As $\psi_t$ is bounded as $t\to 0$, we are justified in taking the 
$\omega$-limit of the function $1/t\to \psi_t(a_0,...,a_p)$, for $t$ 
sufficiently small.

\begin{definition}\label{firstrep} For $p\geq 1$ and $(\A,\HH,\D)$ a $QC^1$ 
$(p,\infty)$-summable spectral triple, set 
\ben \zeta_p(a_0,...,a_p)=
\Omega\mbox{-}\hspace{-.1in}\lim_{1/t\to\infty}\psi_t(a_0,...,a_p),\een
for any (fixed) functional $\Omega$ satisfying the conditions of 
Theorem \ref{niceomega}.
\end{definition}

{\bf Remark} From what we have shown already, $\zeta_p$ is a representative 
of the Hochschild class of the Chern character of $(\A,\HH,\D)$, and
 when $\zeta_p$ is evaluated on a Hochschild cycle, the $\Omega$ limit is a 
 true limit. Thus the value of $\zeta_p$ on Hochschild cycles is independent 
 of the choice of Dixmier functional $\omega$, whether $\omega$ satisfies the 
 extra invariance conditions of Theorem \ref{niceomega} or not.

We need a preliminary result before we can obtain our first formula for 
$\zeta_p$ in terms of the Dixmier trace associated to the functional 
$\Omega$.

\begin{lemma}\label{p1chainrule} Let $(\A,\HH,\D)$ be a $QC^2$ 
$(p,\infty)$-summable spectral triple, with $p\geq 1$  and $\D$ invertible. 
Then
\ben \n[f_p(t\dd),a]-\frac{1}{2}\{f_p'(t\dd),[t\dd,a]\}\n_{(p,1)}
\to 0\ \ \mbox{as}\ \ 
t\to 
0,\een
where $\{T,S\}=TS+ST$.

\medskip
\proof{By Lemma \ref{tracechainrule}, we can estimate the trace norm by
\ben \n[f_p(t\dd),a]-\frac{1}{2}\{f_p'(t\dd),[t\dd,a]\}\n_1\leq Ct^{-p+2}.\een
This completes the proof for $p=1$, and for $p>1$ we 
will employ interpolation as in Lemma \ref{p1bound}.
In Lemma \ref{tracechainrule} we also estimated the operator norm of this 
difference, 
obtaining
\ben \n [f_p(t\dd),a]-\frac{1}{2}\{f_p'(t\dd),[t\dd,a]\}\n_\infty=O(t^2).\een

Applying the interpolation 
inequality
\ben \n T\n_{(p,1)}\leq C\n T\n_1^{1/p}\n T\n^{1-1/p}_\infty,\een
yields
\ben \n[f_p(t\dd),a]-\frac{1}{2}\{f_p'(t\dd),[t\dd,a]\}\n_{(p,1)}\leq 
{\cal C}_{f,p}(t^{-p+2})^{1/p}t^{2-2/p}= 
{\cal C}_{f,p}t\to 0.\een}
\end{lemma}

With these tools in hand, we can now obtain our first Dixmier trace formula 
for 
$\zeta_p$. This result is where we use the invariance properties of the 
Dixmier functional $\Omega$, as this is a necessary condition for 
Theorem \ref{cps2} to hold, at least when $p=1$.

\begin{proposition}\label{usesomega} If $p\geq 1$, $k=\max\{2,p-2\}$ and 
$(\A,\HH,\D)$ is a $QC^k$ 
$(p,\infty)$-summable spectral triple with $\D$ invertible, then for all 
$a_0,...,a_p
\in\A$
\ben \zeta_p(a_0,...,a_p)=p\lambda_p\tau_\Omega(\Gamma 
a_0[F,a_1]\cdots[F,a_{p-1}]\D^{-1}[\dd,a_p]).\een

\medskip
\proof{We begin by noting that we can write $\psi_t(a_0,...,a_p)$ as
\ben  -\lambda_p\tau\left(\Gamma a_0[F,a_1]\cdots[F,a_{p-1}] F
\left(\frac{1}{2}\{f'_p(t\dd),t\delta(a_p)\}+[f_p(t\dd),a_p]-
\frac{1}{2}\{f'_p(t\dd),t\delta(a_p)\}\right)\right).\een
This addition of zero inside the trace is justified as $f'_p(t\dd)$ is trace
 class and $\delta(a_p)$ is bounded. Thus the 
 $\Omega$-limit $\Omega\,\mbox{-}\lim_{1/t\to\infty}\psi_t(a_0,...,a_p)$, is given by 
 the sum of two terms,
\be   -\lambda_p\Omega\mbox{-}\hspace{-.12in}\lim_{1/t\to\infty} \tau\left(\Gamma a_0[F,a_1]\cdots
[F,a_{p-1}] F\frac{1}{2}\{f'_p(t\dd),t\delta(a_p)\}\right)\label{first}\ee
\be-\lambda_p\Omega\mbox{-}\hspace{-.12in}\lim_{1/t\to\infty}\tau\left(\Gamma a_0[F,a_1]\cdots
[F,a_{p-1}]F\left( [f_p(t\dd),a_p]-
\frac{1}{2}\{f'_p(t\dd),t\delta(a_p)\}\right)\right).\label{second}\ee
The second term, (\ref{second}), is zero. To see this, we use the same 
estimates as in Lemma \ref{psitozero}, 
\bean &&\tau\left(\Gamma a_0[F,a_1]\cdots[F,a_{p-1}]F
\left( [f_p(t\dd),a_p]-\frac{1}{2}\{f'_p(t\dd),t\delta(a_p)\}\right)\right)\nno
&\leq &\n\Gamma 
a_0[F,a_1]\cdots[F,a_{p-1}]F\n_{(q,\infty)}\n[f_p(t\dd),a_p]-
\frac{1}{2}\{f'_p(t\dd),t\delta(a_p)\}\n_{(p,1)}\nno
&\leq & {\cal C}_{f_p,p}t\n\Gamma 
a_0[F,a_1]\cdots[F,a_{p-1}]F\n_{(q,\infty)}\to 0\ \ \mbox{as}\ \ t\to 0,\eean
the last inequality following from Lemma \ref{p1chainrule}. Here we again 
replace the $(q,\infty)$ norm by the operator norm when $p=1$.
Hence the  (ordinary) limit of the second term exists and is zero. This means 
that the first term, (\ref{first}), is bounded as $t\to 0$, by 
Lemma \ref{psitozero}, so we have
\bean \zeta_p(a_0,...,a_p) &=& 
-\lambda_p\Omega\mbox{-}\hspace{-.12in}\lim_{1/t\to\infty}\tau(\Gamma 
a_0[F,a_1]\cdots[F,a_{p-1}]F[f_p(t\dd),a_p])\nno
&=&-\lambda_p\Omega\mbox{-}\hspace{-.12in}\lim_{1/t\to\infty}\tau(\Gamma 
a_0[F,a_1]\cdots[F,a_{p-1}]F\frac{1}{2}\{f'_p(t\dd),[t\dd,a_p]\}),\eean
the second line following from Lemma \ref{p1chainrule} and the above argument. 
Using 
\ben f_p'(t\dd)=
-\frac{p}{\Gamma(\frac{p}{2}+1)}t^{p-1}\dd^{p-1}e^{-t^2\D^2},\een
we have
\bean 
\zeta_p(a_0,...,a_p)&=&\frac{p\lambda_p}{2\Gamma(\frac{p}{2}+1)}
\Omega\mbox{-}\hspace{-.12in}\lim_{1/t\to\infty}\tau(\Gamma 
a_0[F,a_1]\cdots[F,a_{p-1}]F\dd^{p-1}t^pe^{-t^2\D^2}[\dd,a_p])\nno
&&\quad+  \frac{p\lambda_p}{2\Gamma(\frac{p}{2}+1)}
\Omega\mbox{-}\hspace{-.12in}\lim_{1/t\to\infty}\tau(\Gamma 
a_0[F,a_1]\cdots[F,a_{p-1}]F[\dd,a_p]\dd^{p-1}t^pe^{-t^2\D^2}).\eean

For $p=1,2$, we use Lemma \ref{bounded} and Theorem \ref{cps2} to obtain
\bean 
\zeta_p(a_0,...,a_p)&=&\frac{p\lambda_p}{2}\tau_\Omega(\Gamma 
a_0[F,a_1]\cdots[F,a_{p-1}]F\dd^{-1}[\dd,a_p])\nno
&&\quad+  \frac{p\lambda_p}{2}\tau_\Omega(\Gamma 
a_0[F,a_1]\cdots[F,a_{p-1}]F[\dd,a_p]\dd^{-1}).\eean
As 
\ben \delta(a_p)\dd^{-1}=
-\dd^{-1}\delta^2(a_p)\dd^{-1}+\dd^{-1}\delta(a_p),\een
 we may commute the $\dd^{-1}$ past $\delta(a_p)$ in the second term, only 
 picking up a term which vanishes under the Dixmier trace. Hence 
\ben 
\zeta_p(a_0,...,a_p)=p\lambda_p\tau_\Omega(\Gamma 
a_0[F,a_1]\cdots[F,a_{p-1}]\D^{-1}[\dd,a_p]).\een

For $p>2$, we use the fact that $\dd^{p-2}e^{-t^2\D^2}$ and 
$\dd e^{-t^2\D^2}$ are trace class, to rewrite $\zeta_p(a_0,...,a_p)$ as
\bean 
&&\frac{p\lambda_p}{2\Gamma(\frac{p}{2}+1)}\Omega\mbox{-}\hspace{-.12in}
\lim_{1/t\to\infty}\tau( t^pe^{-t^2\D^2}\dd[\dd,a_p]\Gamma 
a_0[F,a_1]\cdots[F,a_{p-1}]\dd^{p-2}F)\nno
&&\quad+  \frac{p\lambda_p}{2\Gamma(\frac{p}{2}+1)}
\Omega\mbox{-}\hspace{-.12in}\lim_{1/t\to\infty}\tau(\Gamma 
t^pe^{-t^2\D^2}\dd^{p-2}a_0[F,a_1]\cdots[F,a_{p-1}]F[\dd,a_p]\dd).\eean

Now Lemma \ref{bounded} (and the fact that $\dd$ commutes with $\Gamma$) tells 
us that both of these terms are the
trace of a bounded operator 
times $t^pe^{-t^2\D^2}$. 
So, by Theorem \ref{cps2}, see also \cite{CPS2}, we have
\bean\zeta_p(a_0,...,a_p)&=&\frac{1}{2}p\lambda_p\tau_\Omega(\dd[\dd,a_p]
\Gamma 
a_0[F,a_1]\cdots[F,a_{p-1}]F\dd^{-2})\nno
&&\qquad+\frac{1}{2}p\lambda_p\tau_\Omega(\dd^{-2}\Gamma 
a_0[F,a_1]\cdots[F,a_{p-1}]F[\dd,a_p]\dd).\eean
Since $[F,a_i]$, $[F,\delta(a_i)]$ and $\dd^{-1}$ are in $\LL^{(p,\infty)}$, 
commuting  $\dd$ through these expressions gives, modulo terms of trace class 
which are killed by the Dixmier trace, 
\bean\zeta_p(a_0,...,a_p)&=&\frac{1}{2}p\lambda_p\tau_\Omega([\dd,a_p]\Gamma 
a_0[F,a_1]\cdots[F,a_{p-1}]F\dd^{-1})\nno
&&\qquad+\frac{1}{2}p\lambda_p\tau_\Omega(\dd^{-1}\Gamma 
a_0[F,a_1]\cdots[F,a_{p-1}]F[\dd,a_p]).\eean
In the first term we may cycle $\delta(a_p)$ around to the end using the trace 
property of the Dixmier trace (since $\delta(a_p)$ is bounded while the 
product of the remaining terms is in $\LL^{(1,\infty)}$), while in the second 
we may commute the $\dd^{-1}$ through the product, picking up trace class 
terms from each commutator and these vanish. So
\bea \zeta_p(a_0,...,a_p)&=&\frac{1}{2}p\lambda_p\tau_\Omega(\Gamma 
a_0[F,a_1]\cdots[F,a_{p-1}]\D^{-1}[\dd,a_p])\nno
&+&\frac{1}{2}p\lambda_p\tau_\Omega(\Gamma 
a_0[F,a_1]\cdots[F,a_{p-1}]\D^{-1}[\dd,a_p])\nno
&=&  p\lambda_p\tau_\Omega(\Gamma 
a_0[F,a_1]\cdots[F,a_{p-1}]\D^{-1}[\dd,a_p]).\label{zeta12}\eea}
\end{proposition}

\begin{lemma}\label{cohomega} Let $p\geq 1$, $k=\max\{2,p-2\}$ and let 
$(\A,\HH,\D)$ be a $QC^k$ 
$(p,\infty)$-summable spectral triple with $\D$ invertible. If $\omega$ is any 
Dixmier functional then the multilinear functional $\tilde{\zeta}_p$ defined by
\ben \tilde{\zeta}_p(a_0,...,a_p)=p\lambda_p\tau_\omega(\Gamma a_0[F,a_1]
\cdots[F,a_{p-1}] \D^{-1}[\dd,a_p])\een
represents the Hochschild class of the Chern character.

\medskip

\proof{Let $\Omega$ be a Dixmier functional satisfying the additional 
requirements of Theorem \ref{niceomega}. Then by Proposition \ref{usesomega}, 
the Hochschild class of the Chern character is represented by
\ben \zeta_p(a_0,...,a_p):=p\lambda_p \tau_\Omega(\Gamma a_0[F,a_1]
\cdots[F,a_{p-1}]\D^{-1}[\dd,a_p]).\een
Let $\sum_ia^i_0\otimes a^i_1\otimes\cdots\otimes a^i_p$ be a Hochschild 
cycle, and write the operator
\ben p\lambda_p\sum_i\Gamma a_0[F,a_1]\cdots[F,a_{p-1}]\D^{-1}[\dd,a_p]\een
as a sum $T_1-T_2+iT_3-iT_4$ where $T_i\geq 0$, $i=1,...,4$ and 
$T_i\in\LL^{(1,\infty)}$. Then 
\bean \sum_i\zeta_p(a^i_0,...,a^i_p)&=&
\tau_\Omega(T_1)-\tau_\Omega(T_2)+i\tau_\Omega(T_3)-i\tau_\Omega(T_4)\nno
&=&\Omega\mbox{-}\hspace{-.1in}\;\lim_{t\to\infty}\frac{1}{\log(1+t)}
\int_0^t\mu_t(T_1)dt-\cdots-i\Omega\mbox{-}\hspace{-.1in}\;\lim_{t\to\infty}
\frac{1}{\log(1+t)}\int_0^t\mu_t(T_4)dt\nno
&=&\lim_{t\to\infty}\frac{1}{\log(1+t)}
\int_0^t\mu_t(T_1)dt-\cdots-i\lim_{t\to\infty}\frac{1}{\log(1+t)}
\int_0^t\mu_t(T_4)dt\nno
&=&\omega\mbox{-}\hspace{-.1in}\;\lim_{t\to\infty}\frac{1}{\log(1+t)}
\int_0^t\mu_t(T_1)dt-\cdots-i\omega\mbox{-}\hspace{-.1in}\;\lim_{t\to\infty}
\frac{1}{\log(1+t)}\int_0^t\mu_t(T_4)dt\nno
&=&\sum_i\tilde{\zeta}_p(a^i_0,...,a^i_p).\eean
The equality between the $\Omega$-limit and the true limit follows from 
Lemma \ref{lim} and Proposition \ref{usesomega}, since
\bean \sum_i\zeta_p(a^i_0,...,a^i_p)&=&
\Omega\mbox{-}\hspace{-.12in}\lim_{1/t\to\infty}
\sum_i\psi_t(a^i_0,...,a^i_p)\nno
&=&\lim_{t\to 0}\sum_i\psi_t(a^i_0,...,a^i_p).\eean
Since this is a true limit, any Dixmier functional will also return the same 
value. Note that we are {\em not} asserting that Proposition \ref{usesomega} 
is true for an arbitrary Dixmier functional $\omega$, nor are we asserting 
that $\zeta_p$ and $\tilde{\zeta}_p$ are equal as multilinear functionals. 
What we are asserting is that it makes sense to apply either of $\tau_\omega$ 
or $\tau_\Omega$ to any finite sum of operators of the form
\ben p\lambda_p\sum_i\Gamma a^i_0[F,a^i_1]\cdots[F,a^i_{p-1}]\D^{-1}
[\dd,a^i_p],\ \ a^i_j\in\A,\een
and moreover, that if $c=\sum_ia^i_0\otimes a^i_1\otimes\cdots\otimes a^i_p$ 
is a Hochschild cycle, then $\tau_\omega$ and $\tau_\Omega$ yield the same 
result. The end result of this is that $\tau_\omega-\tau_\Omega$ vanishes on 
all Hochschild cycles. Hence $\tau_\omega$ is cohomologous to $\tau_\Omega$, 
and so $\tau_\omega\in[ICh_F]$.}
\end{lemma}

\begin{corollary}\label{mes} Let $p\geq 1$, $k=\max\{2,p-2\}$ and 
let $(\A,\HH,\D)$ be a $QC^k$ 
$(p,\infty)$-summable spectral triple with $\D$ invertible. For any Hochschild 
cycle $\sum_ia^i_0\otimes a^i_1\otimes \cdots\otimes a^i_p$, the operator 
\ben p\lambda_p\sum_i\Gamma a^i_0[F,a^i_1]\cdots[F,a^i_{p-1}]\D^{-1}
[\dd,a^i_p]\een
is measurable.
\end{corollary}

\subsection{Identification of $[\phi_\omega]$ as the Hochschild Class}

We now come to the cohomological part of the argument where we relate 
$\zeta_p$ to the
functional $\phi_\omega$ appearing in the statement of Theorem \ref{thm8}. 
As mentioned, this part of the proof closely follows \cite[pp477-478]{Va}.

For any choice of Dixmier functional $\omega$, define  
cochains $\zeta_k$, $1\leq k\leq p$, by
\ben \zeta_k(a_0,...,a_p)=
p\lambda_p\tau_\omega(\Gamma 
a_0[F,a_1]\cdots\D^{-1}[\dd,a_k]\cdots[F,a_p]).\een
These are well-defined as the argument of the Dixmier trace in each case is 
an element of $\LL^{(1,\infty)}$ as is readily checked using Lemma 
\ref{sizes}. Note that here we are replacing the definition of $\zeta_p$ given 
in Definition \ref{firstrep}, where we required a Dixmier functional 
satisfying the conditions of Theorem \ref{niceomega}, by the above definition 
using a general Dixmier functional. The two definitions yield cohomologous 
Hochschild cocycles by Lemma \ref{cohomega}.

\begin{lemma}\label{secondlast} Let $p\geq 1$ be integral and suppose that 
$(\A,\HH,\D)$ is a $QC^2$ 
$(p,\infty)$-summable spectral triple with $\D$ invertible.
The cochains $\zeta_1,...,\zeta_p$ are Hochschild cocycles 
which are mutually cohomologous.

\medskip

\proof{We first show that the $\zeta_k$ are Hochschild cocycles. First we need 
to rewrite $\zeta_k$. We wish to rewrite $\D^{-1}\delta(a_k)$ as 
$\delta(a_k)\D^{-1}+$something in $\LL^{(p/2,\infty)}$. First, 
\be [\dd^{-1},T]=-\dd^{-1}[\dd,T]\dd^{-1},\label{commutator}\ee
where  $T=[\dd,a],a,[\D,a]$ or $[F,a]$. So 
\ben \D^{-1}\delta(a_k)=F\delta(a_k)\dd^{-1}-F\dd^{-1}
\delta^2(a_k)\dd^{-1},\een
and the latter term is in $\LL^{(p/2,\infty)}$. Using Lemma \ref{sizes}, we 
see that 
\ben F\delta(a_k)\dd^{-1}=\delta(a_k)\D^{-1}+[F,\delta(a_k)]\dd^{-1}\een
is equal to $\delta(a_k)\D^{-1}$ modulo $\LL^{(p/2,\infty)}$. 
Since each 
$[F,a_j]\in\LL^{(p,\infty)}$, if $T\in\LL^{(p/2,\infty)}$ then we have
\be [F,a_1]\cdots[F,a_{k-1}]T[F,a_{k+1}]\cdots[F,a_p]\in\LL^1.
\label{traceclass}\ee
Hence
\ben\zeta_k(a_0,...,a_p)=p\lambda_p\tau_\omega(\Gamma 
a_0[F,a_1]\cdots[\dd,a_k]\D^{-1}[F,a_{k+1}]\cdots[F,a_p]).\een
To move $\D^{-1}$ all the way to the right, we note that because $F^2=1$, 
$F[F,T]=-[F,T]F$ for all $T\in\NN$, we have
\ben\zeta_k(a_0,...,a_p)=(-1)^{p-k}p\lambda_p\tau_\omega(\Gamma 
a_0[F,a_1]\cdots[\dd,a_k]\dd^{-1}[F,a_{k+1}]\cdots[F,a_p]F).\een
Now 
\ben\dd^{-1}[F,a]=[F,a]\dd^{-1}+[\dd^{-1},[F,a]]= [F,a]\dd^{-1}-\dd^{-1}
[F,\delta(a)]\dd^{-1}\een
and so the operators $\dd^{-1}[F,a]$ and $[F,a]\dd^{-1}$ differ by an element 
of $\LL^{(p/3,\infty)}$ (where for $p<3$ we mean the trace class).

Thus we can move $\D^{-1}$ to the right to obtain
\ben \zeta_k(a_0,...,a_p)=(-1)^{p-k}p\lambda_p\tau_\omega(\Gamma 
a_0[F,a_1]\cdots [\dd,a_k]\cdots [F,a_p]\D^{-1}).\een
Applying Lemma \ref{hochs} and using the trace property of $\tau_\omega$, we 
find that the Hochschild coboundary of 
$\zeta_k$ is given by
\ben (b\zeta_k)(a_0,...,a_{p+1})=
(-1)^{k-1}p\lambda_p\tau_\omega(\Gamma 
a_0[F,a_1]\cdots[\dd,a_k]\cdots[F,a_p][\D^{-1},a_{p+1}]).\een
Repeating the argument of 
Equations \ref{commutator} and \ref{traceclass} shows that this is zero.

The second statement requires that we produce $p$ Hochschild 
$(p-1)$-cocycles  $\eta_k$, 
$k=1,...,p$, such that
\ben b\eta_k(a_0,...,a_p)=\zeta_k-\zeta_{k-1}.\een
The difference on the right hand side is given by 
\ben(-1)^{p-k}p\lambda_p\tau_\omega(\Gamma a_0[F,a_1]\cdots[F,a_{k-2}]
([F,a_{k-1}]\delta(a_k)+\delta(a_{k-1})[F,a_k])[F,a_{k+1}]\cdots[F,a_p]\D^{-
1}).\een
Set $R_{k,k-1}=[F,a_{k-1}]\delta(a_k)+\delta(a_{k-1})[F,a_k]$. Then we have
\be 
[F,\delta(a_{k-1}a_k)]=R_{k,k-1}+a_{k-1}[F,\delta(a_k)]+[F,\delta(a_{k-1})]
a_k.\label{almostder}\ee
So the linear map $a\to [F,\delta(a)]$ is `almost' a derivation. Defining
\ben \eta_k(a_0,...,a_{p-1}):=(-1)^pp\lambda_p\tau_\omega(\Gamma 
a_0[F,a_1]\cdots 
[F,\delta(a_k)]\cdots[F,a_{p-1}]\D^{-1}),\een
it is straightforward to show that $b\eta_k=\zeta_k-\zeta_{k-1}$ using
Equation \ref{almostder} and 
Lemma \ref{hochs}.}
\end{lemma}

\begin{proposition}\label{laststep} Let $p\geq 1$ be integral and suppose that 
$(\A,\HH,\D)$ is a $QC^2$ 
$(p,\infty)$-summable spectral triple with $\D$ invertible.
The cochain 
$\phi_\omega-\frac{1}{p}(\zeta_1+\cdots+\zeta_p)$ is a Hochschild 
coboundary.

\medskip

\proof{We first show that 
\ben \phi_\omega(a_0,...,a_p)=p\lambda_p\tau_\omega(\Gamma a_0[\D,a_1]\cdots
[\D,a_p]\dd^{-p})\een
is equal to the cochain $\tilde{\phi}_\omega$ given by
\ben \tilde{\phi}_\omega(a_0,...,a_p)=p\lambda_p\tau_\omega(\Gamma a_0
[\D,a_1]\dd^{-1}[\D,a_2]\dd^{-1}\cdots[\D,a_p]\dd^{-1}).\een
To do this, we use the argument of Equation (\ref{commutator}) in the last 
Lemma to write
\ben \phi_\omega(a_0,...,a_p)=p\lambda_p\tau_\omega(\Gamma a_0[\D,a_1]\cdots
[\D,a_{p-1}](\dd^{-1}[\D,a_p]\dd^{-p+1}+\dd^{-1}\delta([\D,a_p])\dd^{-p})).\een
The second term is trace class, and so
\ben \phi_\omega(a_0,...,a_p)=p\lambda_p\tau_\omega(\Gamma a_0[\D,a_1]\cdots
[\D,a_{p-1}]\dd^{-1}[\D,a_p]\dd^{-p+1}).\een
Repeating this process of moving one factor of $\dd^{-1}$ to the left at a 
time (which only requires the triple to be $QC^1$) we see that $\phi_\omega=
\tilde{\phi}_\omega$.

Next write 
\ben 
[\D,a_j]\dd^{-1}=[F,a_j]+\delta(a_j)\D^{-1}+[F,\delta(a_j)]\dd^{-1},\een
and observe that by Lemma \ref{sizes} 
$[F,\delta(a_j)]\dd^{-1}\in\LL^{(p/2,\infty)}$. This allows us to replace
\ben [\D,a_j]\dd^{-1}\ \ \mbox{by}\ \ [F,a_j]+\delta(a_j)\D^{-1}\een
 in the formula for $\tilde{\phi}_\omega=\phi_\omega$, using an observation 
 similar to that in Equation (\ref{traceclass}) in the previous Lemma.
Making this substitution will produce $2^p$ functionals, and we will deal 
with them in order of how many terms of the form $\delta(a_j)\D^{-1}$ they 
contain. 
First, we deal with the single functional containing no $\delta(a_j)\D^{-1}$ 
terms,
which is given by
$\lambda_p\tau_\omega(\Gamma a_0[F,a_1]\cdots[F,a_p])$. Now 
$a_0=F[F,a_0]+Fa_0F$, and $F[F,a_0][F,a_1]\cdots[F,a_p]$ is trace class. So
\bean \lambda_p\tau_\omega(\Gamma a_0[F,a_1]\cdots[F,a_p])&=& 
(-1)^p\lambda_p\tau_\omega(\Gamma Fa_0[F,a_1]\cdots[F,a_p]F)\nno
&=& (-1)^{p-1}(-1)^p\tau_\omega(F\Gamma a_0[F,a_1]\cdots[F,a_p]F)\nno
&=& -\lambda_p\tau_\omega(\Gamma a_0[F,a_1]\cdots[F,a_p]).\eean
Hence this functional is zero. The functionals containing precisely one 
$\delta(a_j)\D^{-1}$ term add up to $p^{-1}(\zeta_1+\cdots+\zeta_p)$.

So now we come to the functionals containing two or more terms 
$\delta(a_j)\D^{-1}$. 
So in the following suppose that $\Delta(a)=[F,a]$ or $\delta(a)$, and 
consider a
functional with a total of $l$ terms of the form $\delta(a)\D^{-1}$, 
$2\leq l\leq p$. 
We begin by considering functionals with two consecutive $\delta(a)\D^{-1}$ 
terms. So, modulo
an overall sign arising from moving all powers of $\D^{-1}$ to the right, we 
need to show that
\ben \psi_j(a_0,...,a_p)=\lambda_p\tau_\omega(\Gamma a_0\Delta(a_1)\cdots
\Delta(a_{j-1})
\delta(a_j)\delta(a_{j+1})\Delta(a_{j+2})\cdots\Delta(a_p)\D^{-l})\een
is a coboundary. Now
\ben \delta^2(a_ja_{j+1})=2\delta(a_j)\delta(a_{j+1})+a_j\delta^2(a_{j+1})+
\delta^2(a_j)
a_{j+1},\een
so $\delta^2$ is almost a derivation, and is well-defined on $\A$ since we 
suppose that 
$(\A,\HH,\D)$ is $QC^2$. Setting
\ben \chi_j(a_0,...,a_{p-1})=
\frac{1}{2}(-1)^j\lambda_p\tau_\omega(\Gamma a_0\Delta(a_1)
\cdots\delta^2(a_j)\cdots\Delta(a_{p-1})\D^{-l}),\een
we have, by Lemma \ref{hochs},
\ben (b\chi_j)(a_0,...,a_p)=\psi_j(a_0,...,a_p).\een

So now we are left with functionals in which we do not have two consecutive 
$\delta(a_j)
\D^{-1}$ terms. We will show that such functionals are cohomologous to 
functionals with 
consecutive $\delta(a_j)\D^{-1}$ terms, and so are coboundaries by the 
previous argument.
Again suppose that we have a total of $l$, $2\leq l\leq p$, 
$\delta(a_j)\D^{-1}$ terms.
Consider first 
\ben \xi_{j,j+2}(a_0,...,a_p)=-\lambda_p\tau_\omega(\Gamma a_0\Delta(a_1)\cdots
\Delta(a_{j-1})\delta(a_j)\Delta(a_{j+1})
\delta(a_{j+2})\cdots\Delta(a_p)\D^{-l}),\een
where again $\Delta(a)=[F,a]$ or $\delta(a)$. Let 
\ben \xi_j(a_0,...,a_p)=\lambda_p\tau_\omega(\Gamma a_0\Delta(a_1)\cdots
\Delta(a_{j-1})
\delta(a_j)\delta(a_{j+1})
\Delta(a_{j+2})\cdots\Delta(a_p)\D^{-l}),\een
which is the same as $\xi_{j,j+2}$ except we have swapped the derivations on 
the $j+1$ 
and $j+2$ terms, and introduced an overall minus sign. The difference
$(\xi_{j,j+2}-\xi_j)(a_0,...,a_p)$ is given by
\ben\lambda_p\tau_\omega(\Gamma a_0\Delta(a_1)\cdots\Delta(a_{j+1})\delta(a_j)
(\Delta(a_{j+1})
\delta(a_{j+2})+\delta(a_{j+1})\Delta(a_{j+2}))\cdots\Delta(a_p)\D^{-l}),\een
and this is a coboundary. This is because $\Delta$ and $\delta$ are commuting 
derivations so that
\ben \Delta(\delta(a_{j+1}a_{j+2}))=a_{j+1}\Delta(\delta(a_{j+2}))
+\Delta(\delta(a_{j+1}))a_{j+2}+\Delta(a_{j+1})\delta(a_{j+2})+\delta(a_{j+1})
\Delta(a_{j+2}).\een
Consequently setting
\ben \chi(a_0,...,a_{p-1})=
(-1)^j\lambda_p\tau_\omega(\Gamma a_0\Delta(a_1)\cdots
\Delta(\delta(a_{j+1}))\cdots\Delta(a_{p-1})\D^{-l}),\een
Lemma \ref{hochs} along with the argument following Equation \ref{almostder} 
shows that
\ben (b\chi)(a_0,...,a_p)=\xi_{j,j+2}-\xi_j.\een
Thus any of the functionals containing two or more $\delta(a_j)\D^{-1}$ terms 
are cohomologous to zero. This completes the proof.}
\end{proposition}

This proves Theorem \ref{thm8} for the case where $\D$ has bounded inverse. 
That this is the case
is due to the fact that we can now express the pairing of the Chern character 
with 
Hochschild homology in terms of any 
of the functionals $\zeta_k$, which is the same as employing 
$p^{-1}(\zeta_1+\cdots+\zeta_p)$, and the last Proposition says this is the
same as employing $\phi_\omega$. Theorem \ref{thm8} is also true for the case 
where 
$\D$ does not have bounded inverse; the remaining details are in Appendix 1.

We can also complete the proof of Corollary \ref{ismeas}. Let $\omega$ and 
$\Omega$ be any two Dixmier functionals, and $\sum_ia^i_0\otimes a^i_1\otimes\cdots\otimes a^i_p$ a Hochschild cycle. Then 
\bean \sum_i\phi_\omega(a^i_0,...,a^i_p)&=&
 \sum_i\tau_\omega(\Gamma a^i_0[\D,a^i_1]\cdots[\D,a^i_p]\dd^{-p})\nno
&=&\sum_i\tau_\omega(\Gamma a^i_0[F,a^i_1]\cdots[F,a^i_{p-1}]\D^{-1}[\dd,a_p])\nno
&=&\sum_i\tau_\Omega(\Gamma a^i_0[F,a^i_1]\cdots[F,a^i_{p-1}]\D^{-1}[\dd,a_p])\nno
&=&\sum_i\tau_\Omega(\Gamma a^i_0[\D,a^i_1]\cdots[\D,a^i_p]\dd^{-p})\nno
&=&\sum_i\phi_\Omega(a^i_0,...,a^i_p).\eean
The first equality is the definition of $\phi_\omega$, the second follows from 
Propositions \ref{secondlast} and \ref{laststep}, the third follows from the 
measurability obtained in Lemma \ref{cohomega} and Corollary \ref{mes}, and 
the final two equalities follow from Propositions \ref{secondlast} and 
\ref{laststep} and the definition. Hence the operator
\ben \sum_i\Gamma a^i_0[\D,a^i_1]\cdots[\D,a^i_p]\dd^{-p}\een
is measurable, and Corollary \ref{ismeas} is proved.

\section{Appendix}\label{appone}

Our chief remaining task is to determine the effects on our 
representative $\phi_\omega$ of the Hochschild class of the Chern character 
of replacing $(\A,\HH,\D)$ by $(\A,\HH^2,\D_m)$.

In \cite[III.1.$\beta$,Proposition 15]{BRB}, Connes shows that if 
$\phi\in Z^n_\lambda(\A)$ is a cyclic cocycle, then $S\phi:=\phi\#\s$ is a 
Hochschild 
coboundary. Here $\#$ is the cup product, \cite[pp 191-193]{BRB}, and $\s$ 
defined by $\s(1,1,1)=1$ is the 
cyclic cocycle generating the cyclic cohomology of ${\C}$. It is important to 
realise 
that $\s$ is a Hochschild coboundary.

To define the periodicity operator on arbitrary cyclic cochains, one must 
introduce antisymmetrisation and some normalisation constants. This is not an 
appropriate procedure for Hochschild cochains, and it is in fact simply the 
cup product by the cyclic cocycle (Hochschild {\em coboundary}) $\s$ which is 
important for us. Consequently, for any Hochschild cycle $\phi$, we shall 
denote by $S\phi$ the Hochschild cocycle $\phi\#\s$. Note that this is not the 
usual definition 
of the periodicity operator $S$, but our definition coincides with the usual 
definition 
on cyclic cocycles.
The important point is that if $\phi$ is a Hochschild cocycle, 
then $\phi\#\s$ is a Hochschild coboundary, \cite[p 194]{BRB}.

Our strategy is to show that the representative of the Hochschild class of 
the Chern 
character we obtain in Theorem \ref{thm8} when we use the operator 
\ben \bma \D & m\\ m & -\D\ema\een
differs from our stated result by Hochschild coboundaries.

\begin{definition}\label{functionals} Let $p\geq 1$ be integral and suppose 
that 
$(\A,\HH,\D)$ is a $QC^2$ $(p,\infty)$-summable spectral triple. Let 
$(\A,\HH^2,\D_m)$ 
be the `double' of $(\A,\HH,\D)$, and define
\bean \phi_\omega(a_0,...,a_p) &=&
 \lambda_p\tau_\omega(\Gamma a_0[\D,a_1]\cdots 
[\D,a_p](m^2+\D^2)^{-p/2})\nno
&=& \lambda_p\tau_\omega(\Gamma a_0[\D,a_1]\cdots [\D,a_p](1+\D^2)^{-p/2})\eean
\ben \phi_{\omega m}(a_0,...,a_p)=
\lambda_p\tau_\omega(\Gamma a_0[\D_m,a_1]\cdots 
[\D_m,a_p]|\D_m|^{-p})\een
\bean \tilde{\phi}^k_\omega(a_0,...,a_k)&=&
\lambda_p\tau_\omega(\Gamma a_0[\D,a_1] 
\cdots[\D,a_k](m^2+\D^2)^{-p/2})\nno
&=& \lambda_p\tau_\omega(\Gamma a_0[\D,a_1] \cdots[\D,a_k]
(1+\D^2)^{-p/2}).\eean
\end{definition}

The equalities in the definition follow from 
\ben(m^2+\D^2)^{-1}-(n^2+\D^2)^{-1}=
(n^2-m^2)(m^2+\D^2)^{-1}(n^2+\D^2)^{-1},\een
which, by the BKS inequality \cite{BKS}, implies that
\ben (m^2+\D^2)^{-p/2}-(n^2+\D^2)^{-p/2}\in\LL^1.\een 
Hence $\tau_\omega(A(m^2+\D^2)^{-p/2})=\tau_\omega(A(n^2+\D^2)^{-p/2})$ for 
all 
bounded $A\in\NN$ and $n,m>0$.

It is straightforward to show using Lemma \ref{hochs} and/or Lemma 
\ref{iscocycle} that all of the functionals in Definition \ref{functionals}
are Hochschild cocycles. The explicit formula for $\phi\#\s$ where $\phi$ is 
any  of the above ($n$-)cocycles, is \cite[p 193]{BRB},
\bean (\phi\#\s)(a_0,...,a_{n+2})&=& \phi(a_0a_1a_2da_3\cdots da^{n+2}T)\nno
&&+\phi(a_0da_1 (a_2a_3)da_4\cdots da_{n+2}T)+\cdots\nno
&&+\phi(a_0da_1\cdots da_{i-1}(a_ia_{i+1})da_{i+2}\cdots da_{n+2}T)+\cdots\nno
&&+\phi(a_0da_1\cdots da_n(a_{n+1}a_{n+2})T),\eean
where $da$ denotes $[\D,a]$ and we have written $T$ generically for 
$(1+\D^2)^{-p/2}$ or $|\D_m|^{-p}$ etc.

We can now state the main result of the Appendix.

\begin{proposition}\label{nounit} Let $p\geq 1$ be integral and suppose that 
$(\A,\HH,\D)$ is a $QC^2$ $(p,\infty)$-summable spectral triple. Let 
$(\A,\HH^2,\D_m)$ 
be the `double' of $(\A,\HH,\D)$. Then for all $a_0,...,a_p\in\A$
\bean \phi_{\omega m}(a_0,...,a_p)&=&\phi_\omega(a_0,...,a_p)\nno
&+&\sum_{i=1}^{[p/2]}(-1)^im^{2i}\frac{1}{i!}(S^i\tilde{\phi}^{p-2i}_\omega)
(a_0,...,a_p).\eean

\medskip

\proof{We begin by defining a collection of operators $\hat{S}^i$, $i\geq 1$, 
which we 
will use to work with elements of $\Omega^*_\D(\A)$, the graded algebra 
generated by $\A$ and $[\D,\A]$, rather than with the cocycles. We 
define $\hat{S}:\A^{\otimes n+1}\to\Omega^{n-2}_\D(\A)$, for any $n$ by
\ben \hat{S}(a_0)=\hat{S}(a_0,a_1)=0,\een
\ben \hat{S}(a_0,...,a_n)=
\sum_{i=1}^{n-1}a_0d(a_1)\cdots d(a_{i-1})a_ia_{i+1} 
d(a_{i+2})\cdots d(a_n).\een
Here and below we write $d(a)=[\D,a]$. To define `powers' of $\hat{S}$, we 
employ the 
inductive definition
\bean \hat{S}^k(a_0,...,a_n)&=&
\hat{S}^{k-1}(\sum_{i=1}^{n-1}a_0d(a_1)\cdots d(a_{i-1})
a_ia_{i+1} d(a_{i+2})\cdots d(a_n))\nno
&=&\hat{S}^{k-2}(\sum_{i=1}^{n-1}\hat{S}(a_0,...,a_{i-1})a_ia_{i+1}d(a_{i+2})
\cdots d(a_n))\nno
&+&\hat{S}^{k-2}(\sum_{i=1}^{n-1}a_0d(a_1)\cdots d(a_{i-1})\hat{S}(a_ia_{i+1},
a_{k+2},...,a_n))\nno
&=&\sum_{j=0}^{k-1}\bca k-1\\ j\eca\sum_{i=1}^{n-1}\hat{S}^j(a_0,...,a_{i-1})
\hat{S}^{k-j-1}(a_ia_{i+1},...,a_n).\eean

It is tedious but not difficult to check that 
\be (S^i\phi)(a_0,...,a_n)=\phi(\hat{S}^i(a_0,...,a_n)),\label{tedious}\ee
for any of the Hochschild cocycles defined in Definition \ref{functionals} 
(regarded as  functionals on $\Omega^*_\D(\A)$).

We claim that for any $n\geq 0$ the product $a_0[\D_m,a_1]\cdots[\D_m,a_n]$ is 
given by
\be \bma \begin{array}{l}a_0d(a_1)\cdots d(a_n)\\+\sum_{i=1}^{[n/2]}
\frac{1}{i!}
(-1)^im^{2i}\hat{S}^i(a_0,...,a_n)\end{array} & \begin{array}{l} ma_0d(a_1)
\cdots 
d(a_{n-1})a_n\\+\sum_{i=1}^{[(n-1)/2]}m^{2i+1}(-1)^i\frac{1}{i!}\hat{S}^i
(a_0,...,a_{n-1})a_n\end{array}\\ 0 & 0\ema.\label{claim}\ee
Indeed, this is easy to verify for $n=1,2$. So if we suppose it to be true 
for all 
$k<n$ then using
\ben [\D_m,a_n]=\bma d(a_n) & ma_n\\ -ma_n & 0\ema,\een
we find that $a_0[\D_m,a_1]\cdots[\D_m,a_n]$ is given by 
(writing $c_i=\frac{1}{i!}(-1)^im^{2i}$)

\bean &&\!\bma \begin{array}{l}a_0d(a_1)\cdots d(a_{n-1})\\
+\sum_{i=1}^{[(n-1)/2]}c_i\hat{S}^i(a_0,...,a_{n-1})\end{array} & 
\begin{array}{l} 
ma_0d(a_1)\cdots d(a_{n-2})a_{n-1}\\+m\sum_{i=1}^{[(n-2)/2]}
c_i\hat{S}^i(a_0,...,a_{n-2})a_{n-1}
\end{array}\\ 0 & 0\ema [\D_m,a_n]\nno
&=&\bma \begin{array}{l} a_0d(a_1)\cdots d(a_n)\\ +\sum_{i=1}^{[(n-1)/2]}
c_i\hat{S}^i(a_0,...,a_{n-1})d(a_n)\\
-m^2\sum_{i=0}^{[(n-2)/2]}c_i\hat{S}^i(a_0,...,a_{n-2})a_{n-1}a_n\end{array} &
\begin{array}{l} ma_0d(a_1)\cdots d(a_{n-1})a_n\\
+\sum_{i=1}^{[(n-1)/2]}c_i\hat{S}^i(a_0,...,a_{n-1})a_n\end{array}
\\ 0 & 0\ema.\eean
In order to simplify this expression we note that 
\be \hat{S}(a_0,...,a_{n-1})d(a_n)=
\hat{S}(a_0,...,a_n)-a_0d(a_1)\cdots d(a_{n-2})
a_{n-1}a_n,\label{one}\ee
and for $i>1$
\ben \hat{S}^i(a_0,...,a_{n-1})d(a_n)=\hat{S}^i(a_0,...,a_n)-i\hat{S}^{i-1}
(a_0,...,a_{n-2})a_{n-1}a_n.\een

To see this, one first verifies the statement for $i=2$ (which is 
straightforward using
Equation \ref{one} and a calculation similar to that below), and then we use 
induction. 
The computation is as follows.

\bea \hat{S}^{k+1}(a_0,...,a_n) & = &
 \sum_{j=0}^k\bca k \\ j\eca\sum_{i=1}^{n-1}
\hat{S}^j(a_0,...,a_{i-1})\hat{S}^{k-j}(a_ia_{i+1},...,a_n)\nno
&=&\sum_{j=0}^{k-1}\bca k\\ j\eca
\sum_{i=1}^{n-1}\hat{S}^j(a_0,...,a_{i-1})\hat{S}^{k-j}
(a_ia_{i+1},...,a_{n-1})d(a_n)
\nno
&+&\sum_{j=0}^{k-1}\bca k\\ j\eca(k-j)
\sum_{i=1}^{n-1}\hat{S}^j(a_0,...,a_{i-1})
\hat{S}^{k-j-1}(a_ia_{i+1},...,a_{n-2}) a_{n-1}a_n\nno
&+&\sum_{i=1}^{n-2}\hat{S}^k(a_0,...,a_{i-1})a_ia_{i+1}d(a_{i+2})\cdots d(a_n) 
+\hat{S}^k(a_0,...,a_{n-2})a_{n-1}a_n\label{induct}\\
&=& \sum_{j=0}^{k}\bca k\\ j\eca
\sum_{i=1}^{n-2}\hat{S}^j(a_0,...,a_{i-1})\hat{S}^{k-j}
(a_ia_{i+1},...,a_{n-1})d(a_n) \nno
&+& k\sum_{j=0}^{k-1}\bca k-1\\ j\eca\sum_{i=1}^{n-1}\hat{S}^j(a_0,...,a_{i-1})
\hat{S}^{k-j-1}(a_ia_{i+1},...,a_{n-2}) a_{n-1}a_n\nno
&+& \hat{S}^k(a_0,...,a_{n-2})a_{n-1}a_n\label{zero}\\
&=& \hat{S}^{k+1}(a_0,...,a_{n-1})d(a_n)+(k+1)
\hat{S}^{k}(a_0,...,a_{n-2})a_{n-1}a_n.
\label{result}\eea
The first line here follows from the definition. In \ref{induct} we apply the 
inductive 
hypothesis to the second term in each product, for $j\neq k$, and for $j=k$ we 
split 
the sum into the first $n-2$ terms, and the $(n-1)$-st. For $j\neq k$ we 
notice that 
the $(n-1)$-st term of the sum is zero, by the definition of $\hat{S}$, so 
in \ref{zero}
we collect all these sums of $n-2$ terms. We also use the combinatorial 
identity
\ben (k-j)\bca k\\ j\eca=k\bca k-1\\ j\eca.\een
Finally, applying the definition of $\hat{S}$ we obtain the result 
\ref{result}.

Thus we have
\bean &&\sum_{i=1}^{[(n-1)/2]}c_i\hat{S}^i(a_0,...,a_{n-1})d(a_n)
-m^2\sum_{i=0}^{[(n-2)/2]}c_i\hat{S}^i(a_0,...,a_{n-2})a_{n-1}a_n\nno
&&\nno
&=&\left\{\begin{array}{ll}\sum_{i=1}^{[n/2]}(c_i\hat{S}^i(a_0,...,a_{n-1})
d(a_n)+ic_i\hat{S}^{i-1}(a_0,...,a_{n-2})a_{n-1}a_n) & n\ \mbox{odd}\\
& \\

\begin{array}{l}\sum_{i=1}^{[(n-1)/2]}(c_i\hat{S}^i(a_0,...,a_{n-1})d(a_n)+
ic_i\hat{S}^{i-1}(a_0,...,a_{n-2})a_{n-1}a_n\\ +m^n(-1)^{n/2}
\frac{1}{([n/2]-1)!}
\hat{S}^{[n/2]-1}(a_0,...,a_{n-2})a_{n-1}a_n\end{array} & n\ \mbox{even}
\end{array}
\right.\nno
&=&\left\{\begin{array}{ll}\sum_{i=1}^{[n/2]}c_i\hat{S}^i(a_0,...,a_n) & n\ 
\mbox{odd}\\  & \\
\begin{array}{l}\sum_{i=1}^{[(n-1)/2]}c_i\hat{S}^i(a_0,...,a_n)\\
+m^n(-1)^{n/2}\frac{1}{([n/2]-1)!}\hat{S}^{[n/2]-1}(a_0,...,a_{n-2})a_{n-1}a_n
\end{array} & n\ \mbox{even}\end{array}\right..\eean
In the odd case we have used $[(n-1)/2]=[n/2]$, and we are left with the even 
case. 
For this we note that

\ben k\hat{S}^{k-1}(a_0,...,a_{2k-2})a_{2k-1}a_{2k}=
\hat{S}^k(a_0,...,a_{2k}),\een
since $\hat{S}^k(a_0,...,a_{2k-1})=0$, so 
\ben m^n(-1)^{n/2}
\frac{1}{([n/2]-1)!}\hat{S}^{[n/2]-1}(a_0,...,a_{n-2})a_{n-1}a_n=
m^{2[n/2]}(-1)^{[n/2]}\frac{1}{[n/2]!}\hat{S}^{[n/2]}(a_0,...,a_n).\een

This completes the inductive step and proves the claim \ref{claim}.
Putting \ref{claim} together with \ref{tedious} now completes the proof.}
\end{proposition}

Thus the Hochschild class of the Chern character can be represented by the 
cocycle 
\ben \phi_\omega(a_0,...,a_p)=\lambda_p\tau_\omega(\Gamma a_0[\D,a_1]\cdots 
[\D,a_p](1+\D^2)^{-p/2}),\een
the other contributions appearing in Proposition \ref{nounit} all being 
coboundaries with no effect on the Hochschild class.

\end{document}